\documentclass{article}

\usepackage{PRIMEarxiv}

\usepackage[utf8]{inputenc} 
\usepackage[T1]{fontenc}    
\usepackage{hyperref}       
\usepackage{url}            
\usepackage{booktabs}       
\usepackage{amsfonts}       
\usepackage{nicefrac}       
\usepackage{microtype}      
\usepackage{lipsum}
\usepackage{fancyhdr}       
\usepackage{graphicx}       
\graphicspath{{media/}}     

\usepackage{amsmath}
\usepackage{siunitx}
\usepackage{longtable,tabularx}
\usepackage{enumitem}
\usepackage{algorithm}
\usepackage{algpseudocode}
\usepackage{graphicx,subcaption}

\title{Exploiting Scaling Constants to Facilitate the Convergence of Indirect Trajectory Optimization Methods
}

\author{
  Minduli C. Wijayatunga \\
  Phd Candidate\\
  Te Pūnaha Ātea - Space Institute \\
  The University of Auckland \\
  Auckland\\
  \texttt{mwij516@aucklanduni.ac.nz} \\
   \And
  Roberto Armellin\\
  Professor \\
  Te Pūnaha Ātea - Space Institute \\
  The University of Auckland \\
  Auckland\\
  \texttt{roberto.armellin@auckland.ac.nz} \\
     \And
  Laura Pirovano\\
  Research Fellow \\
  Te Pūnaha Ātea - Space Institute \\
  The University of Auckland \\
  Auckland\\
  \texttt{laura.pirovano@auckland.ac.nz } \\}
  
\begin{document}
\maketitle

\begin{abstract}
  This note develops easily applicable techniques that improve the convergence and reduce the computational time of indirect low thrust trajectory optimization when solving fuel- and time-optimal problems.
For solving fuel optimal (FO) problems, a positive scaling factor-$\Gamma_{TR}$ -is introduced based on the energy optimal (EO) solution to establish a convenient profile for the switching function of the FO problem. This negates the need for random guesses to initialize the indirect optimization process. Similarly, another scaling factor-$\beta$-, is introduced when solving the time-optimal (TO) problem to connect the EO problem to the TO.   The developed methodology for the TO problem was crucial for the GTOC11 competition. Case studies are conducted to validate the solution process in both TO and FO problems. For geocentric cases, the effect of eclipses and $J_2$ perturbations were also considered. The examples show that EO can provide a good guess for TO and FO problems and that introducing the constants can reduce the initial residuals and improve convergence. It is also shown that the equation for the Lagrangian multiplier of mass and the associated boundary condition can be ignored for both FO and TO cases without affecting optimality. This simplification reduces the problem dimensions and improves efficiency.
\end{abstract}


\section{Introduction}
Low thrust propulsion systems have gained traction in recent years as they are capable of significantly reducing the propellant consumption of a mission by providing a higher specific impulse \cite{1}. However, optimizing continuous  {low-thrust} trajectories is significantly more complicated than optimizing trajectories involving chemical propulsion. When optimizing a chemical propulsion trajectory, only a finite number of variables need to be considered (i.e., number, magnitude, and direction of impulses). In contrast,  {low-thrust} optimization requires a continuous steering law to be determined for the entire trajectory while satisfying system constraints and operational restrictions \cite{2}. Furthermore, the nonlinear and non-convex nature of  {low-thrust} dynamics, environmental perturbations, and the presence of local minima introduce further complications to the optimization process \cite{3}. \par 
A  {low-thrust} trajectory optimization problem can be set up to maximize or minimize a particular cost function (i.e. minimizing time of flight \cite{Geller2017,Taheri2021,Wu20212} or minimizing the fuel consumption \cite{Altitude-dependent2010c,Wu2021,Rasotto2015,Kechichian1997,Woollands2020}). Many numerical and analytical approaches have been introduced to solve the  {low-thrust} trajectory optimization over the years, based on classical optimal control \cite{Taheri2021, Jin2020, Ravikumar2020, Armellin2018c, Altitude-dependent2010c, Wu2021, Russell2007a, Wall2009}. Traditionally, these methods are categorized into direct and indirect methods. \par 
For direct methods, the optimal control problem (OCP) is transcribed into a nonlinear programming problem where the objective function is directly optimized using Karush-Kuhn-Tucker (KKT) optimality conditions \cite{3}. In indirect methods, the OCP is converted into a  multi-point boundary value problem (MPBVP)  and solved by applying the Pontryagin minimum/maximum principle (PMP) \cite{prussing_2010}. Direct methods are simpler to initialize and have more robust convergence. However, the optimality of obtained solution depends on how the control is parameterized. In contrast, indirect methods ensure the satisfaction of the first-order optimality conditions. Furthermore, they can provide theoretical knowledge on the physical and mathematical features of the problem \cite{3}. However, the indirect optimization method also faces significant challenges. These include the need for an initial guess of the costates, convergence issues due to discontinuities, and difficulties handling the state constraints \cite{Taheri2021}.\par Methods have been developed to calculate the initial guess for indirect methods analytically, solving the problem using warm start methods \cite{Wu2021, Kim2008}. Randomly guessing the initial costates on an eight-dimensional unit hyper-sphere (adjoint normalization method) and generating an initial guess within a specified range have also been utilized to solve this issue \cite{Altitude-dependent2010c, Taheri2018, Jiang2012}.  \par 
Continuation methods have been frequently utilized to overcome the difficulties caused by discontinuities and state constraints \cite{Bertrand, Haberkorn, Rasotto2015, Taheri2018}. Continuation involves solving a series of subproblems that become incrementally more complex, ultimately leading to the desired solution \cite{mansell2018adaptive}. 
For solving the fuel-optimal (FO) problem, a more straightforward subproblem has been constructed with the objective function of the optimization taken to be the square of the thrust magnitude \cite{Bertrand, Martinon2007}. This problem, known as the energy-optimal (EO) problem, has been used many times to solve FO trajectories since its introduction \cite{Wu2021, Wu20212, Taheri2018,tang2018fuel}. Several ways of connecting the FO and EO problems have been studied in detail thus far \cite{Rasotto2015, ZHU201798, Haberkorn, Wu2021, Jiang2017, tang2018fuel, Altitude-dependent2010c}. However, methods of connecting the time-optimal (TO) and EO problems remain relatively unexplored.
\par 
In recent years, there has been some exploration into the effect of introducing positive scaling factors to the objective function to improve various aspects of the problem formulation \cite{Jiang2012, Ayyanathan2022}. In \cite{Jiang2012}, a factor is introduced in the FO problem to restrict the search space for all Lagrangian multipliers in $[-1,1]$. This approach is particularly beneficial when global optimization is used to find a first-guess solution. In \cite{Ayyanathan2022}, the scaling factor is used to convert the costates between coordinate systems. However, using scaling factors to facilitate the continuation method and connect the desired problem to the EO problem has not been explored thus far. \par 
This note develops novel and easily applicable techniques that, when applied with smoothing and continuation, eliminate the need to generate random guesses for the initial Lagrangian multipliers. Consequently, our approach can improve convergence and reduce the computational time of indirect methods, even in the presence of perturbations and eclipses. \par 
For solving the FO problem, a positive scaling factor, referred to as the thrust threshold ($\Gamma_{TR}$), is introduced. This parameter is calculated through the EO solution and used to establish a convenient profile for the switching function of the FO problem. This approach has been recently introduced to optimize collision avoidance maneuvers when the variation of the spacecraft mass is negligible \cite{palermoCAM}. Here we extend it to the general case, including mass variation by using the total $\Delta v$ as objective and control in acceleration so that the dynamics of the Lagrangian multiplier of the mass can be neglected. \par 
A similar strategy is used to solve the TO problem. The EO problem is used to estimate the optimal time of flight and to provide an initial guess for the Lagrangian multipliers. Moreover, a positive scaling factor ($\beta_t$) is introduced in the TO objective function to connect it to the EO problem. This constant is set such that the constraint on the final Hamiltonian is satisfied by the first guess solution. Additionally, we show that even in the TO problem, the costate of the mass and its boundary condition can be dropped, thus further simplifying the problem.  {The developed methodology for this problem was crucial during the GTOC11 competition in which our team \lq \lq theAntipodes" ranked third \footnote{Results of the competition are available at \url{https://gtoc11.nudt.edu.cn/}}}\par 
The rest of this note is organized as follows. Section II discusses the spacecraft dynamics associated with the optimization problems. The problem formulation for the FO case is presented in Section III. In section IV, the TO problem is formulated. Section V reviews the results obtained by solving the FO and TO problems for geocentric and heliocentric trajectories. Lastly, concluding remarks are given in Section VI.

\section{Spacecraft Dynamics}\label{dynamics}

The motion of a spacecraft can be modeled as
\begin{equation}\label{d1}
\boldsymbol{\dot{x}}=\boldsymbol{A}\left(\boldsymbol{x} \right) + \boldsymbol{B}\left(\boldsymbol{x} \right)  \frac{T_{max}}{m_0} \Gamma  \widehat{\boldsymbol{\alpha}}, 
\end{equation}
where the control acceleration magnitude is given by $\frac{T_{max}}{m_0} \Gamma $ where $ \Gamma \in \left[ 0, \frac{m_0}{m(t)} \right] $. Its direction is given by $\widehat{\boldsymbol{\alpha}}$. The propulsion system is assumed to have a maximum thrust magnitude of $T_{max}$ and a specific impulse of $I_{sp}$. $m_0$ denotes the initial mass of the spacecraft. The expressions for matrices $\boldsymbol{A}$ and $\boldsymbol{B}$ are dependent on the coordinate system used. 
 {The modified equinoctial elements are used throughout the note, hence the state vector is defined as $\boldsymbol{x} = [p, f, g, h, k, L]$. The definition of each coordinate of $\boldsymbol{x}$ and $\boldsymbol{A}$ and $\boldsymbol{B}$ matrices are provided in \cite{Betts2000VeryLT}. }
 \par 
 The mass of the spacecraft is kept track of by noting that   {according to the rocket equation \cite{prussing_2010},
 \begin{equation}\label{m1}
     m(t) = m_0 \exp\left(-{\frac{\Delta v(t)}{I_{sp}g_0}}\right), \ \text{ where} \   \Delta v(t) =  \frac{T_{max}}{m_0} \int^{t}_{t_0}  \Gamma(\xi) \  d\xi.
 \end{equation}}
and $g_0$ denotes the standard acceleration of gravity at sea level  {(9.80665 m/s$^2$)}. \par 
 The effect of eclipses and $J_2$ perturbations are introduced to the spacecraft dynamics for geocentric transfers. The effect of $J_2$ in the radial, transverse and normal (RTN) coordinates can be defined as follows using modified equinoctial elements, where $J_2 = 1.08262668 \times 10^{-3}$  {, $h,k$ and $L$ are the modified equinoctial elements in $\boldsymbol{x}$ and $\mu$ is the standard gravitational parameter of Earth ($ \SI{3.986e5}{\kilo\meter\cubed\per\second\squared}$). $R_E$ denotes the radius of Earth and is taken to be $6378.1370$ km. $r$ depicts the distance from the spacecraft to the centre of the Earth. }
  
  \begin{equation}
   \Delta \boldsymbol{J}_2 = \frac{-3\mu  J_2 R_E^2}{2 r^4(1+ h^2 + k^2)^2} \begin{bmatrix}  1 - 12 (h\sin{L} - k\cos{L} )^2 \\ 8(h\sin{L} - k\cos{L})(h\cos{L} + ksin{L})\\ 4(h\sin{L} - k\cos{L})(1 - h^2 - k^2) \end{bmatrix} 
\end{equation}
A spacecraft in a geocentric orbit experiences eclipses when the shadow of the Earth obstructs the Sun. This causes a momentary loss of thrust. The scale factor $\nu$ is introduced to represent eclipses, where a full eclipse is denoted by $\nu = 1$, a partial eclipse by $0 < \nu <1$. $\nu =0$ when no eclipse is present.


There are various methods of determining $\nu$ in the literature \cite{don, Vallado2013, Jonathan2019}. In this note, it is determined by the eclipse smoothing technique described in \cite{Jonathan2019}. Hence, for geocentric transfers that take eclipses and the effect of $J_2$ into consideration, the spacecraft dynamics are modeled as 

\begin{equation}\label{d2}
\boldsymbol{\dot{x}}=\boldsymbol{A}\left(\boldsymbol{x} \right)+ \boldsymbol{B}\left(\boldsymbol{x}\right)  \frac{T_{max}}{m_0} \left(1-\nu \right) \Gamma   \widehat{\boldsymbol{\alpha}} + \boldsymbol{B}\left(\boldsymbol{x}\right) \Delta \boldsymbol{J}_2  
\end{equation}
where  the mass variation is calculated using Eq. \eqref{m1}. 
\section{Fuel-Optimal Problem }

The FO problem entails minimizing 
\begin{equation}
    J = \Gamma_{TR} \frac{T_{max}}{m_0} \int^{t_1}_{t_0}   \Gamma_f(t) \  dt, 
    \label{lambdanonzerofopt2}
\end{equation}

subject to the dynamics given in Eq. \eqref{d1}. $t_0$ and $t_1$ denote the start and end times of the mission, respectively. The subscript $f$ is used to indicate quantities in the FO problem. $\Gamma_{TR}$ is a positive constant whose value will be determined later.
The terminal constraints that must hold are
\begin{equation}\label{conditions2}
\boldsymbol{x}(t_0) = \boldsymbol{x_0}, \  \boldsymbol{x}(t_1) = \boldsymbol{x_1}.
\end{equation} 

where $\boldsymbol{x_0}$ is the starting position of the spacecraft and $\boldsymbol{x_1}$ denotes its target destination.
This problem can be transformed into a two-point boundary value problem (TPBVP) using Pontryagin's minimum/maximum principle. 
First, the Hamiltonian of this problem is built as follows, by introducing the costate vector $\boldsymbol{\lambda_x} = [\lambda_p, \lambda_f, \lambda_g, \lambda_h, \lambda_k, \lambda_L]^T$
\begin{equation}\label{H1}
     \mathcal{H} =  \boldsymbol{\lambda_x}^T\left[ \boldsymbol{A} + \boldsymbol{B}  \frac{T_{max}}{m_0}\Gamma_f \widehat{\boldsymbol{\alpha}}_f\right] +  \Gamma_{TR}\frac{T_{max}}{m_0} \Gamma_f .
\end{equation}
The problem dynamics and the costate differential equations are then be written as  {$\boldsymbol{\dot{x}} = \frac{\partial{\mathcal{H}}}{\partial{ \boldsymbol{\lambda}}} $ and $\boldsymbol{\dot{\lambda}}_{\boldsymbol{x}} = -\frac{\partial{\mathcal{H}}}{\partial{ \boldsymbol{x}}}$.} The optimal thrust direction that minimizes the Hamiltonian given in Eq. \eqref{H1} is 
\begin{equation}
    \widehat{\boldsymbol{\alpha}}_f^* = -\frac{ \boldsymbol{B^T \lambda_x}}{\|\boldsymbol{B^T \lambda_x}\|}
\end{equation}
Substituting $\widehat{ \boldsymbol{\alpha}}_f^*$ to Eq. \eqref{H1}, the 
$\Gamma_f$ that minimizes the Hamiltonian can be determined as: 
\begin{equation}\label{foptthrust}
\begin{aligned}
\text{If}\ \Gamma_{TR} -\|\boldsymbol{B^T \lambda_x}\|> 0 &\Rightarrow \Gamma_f^* = 0   \\
\text{If}\ \Gamma_{TR} -\|\boldsymbol{B^T \lambda_x}\|< 0 &\Rightarrow \Gamma_f^* = \frac{m_0}{m(t)}.   \\
     \end{aligned}
\end{equation}
Hence $\Gamma_f^*$ can be defined as: 
\begin{equation}
    \Gamma_f^* = \frac{m_0}{2m(t)}  \left[1- \text{sgn}(\rho)\right], \ 
\end{equation}
where the switching function ($\rho$) is defined as 
\begin{equation}\label{FOSF}
  \rho = \Gamma_{TR} -\|\boldsymbol{B^T \lambda_x}\|.
\end{equation}
Smoothing is implemented on $\Gamma_f^*$ to avoid numerical difficulties associated with the discontinuous behavior of the thrust \cite{Rasotto2015, Ayyanathan2022, Lizia2014}.  One such smoothed representation of $\Gamma_f^*$ is
\begin{equation}
    \Gamma_f^* =   \frac{m_0}{2m(t)} - \frac{m_0}{2m(t)} \tanh \left({\frac{\rho}{1- k}}\right),
\end{equation}
where the smoothing parameter $k \in [0,1)$. By increasing $k$ from 0 to 1, the thrust profile approximates more accurately the optimal discontinuous profile.  When this smoothed representation of thrust is used in FO problems, it shall be known as the Smoothed Fuel-Optimal (SFO) problem from now on.
The initial ($\boldsymbol{x(t_0)}$) and final states ($\boldsymbol{x(t_1)}$) of the system must satisfy Eq. \eqref{conditions2}. Thus, the problem now consists of finding the initial lagrangian multipliers $\boldsymbol{\lambda_{x}}(t_0)$ that satisfy
\begin{equation}\label{Phi1}
    \boldsymbol{\Phi}(\boldsymbol{\lambda_{x}}(t_0)) = [ \boldsymbol{x}(t_1) - \boldsymbol{x_1}] = \boldsymbol{0},
\end{equation}
which is known as a shooting function. 
Eq. \eqref{Phi1} can then be solved using a nonlinear solver, such as the multidimensional root finding method in the GNU Scientific Library \cite{galassi2002gnu}.

\subsection{Energy-Optimal Problem } \label{eoptsection}

A non-intuitive initial guess is required for the convergence of the FO solution when it is solved via nonlinear methods. A typical approach to obtain a sufficiently accurate initial guess is to solve the EO problem \cite{Altitude-dependent2010c}. In this note, this approach is further enhanced by exploiting $\Gamma_{TR}$ to link the EO solution to the FO problem optimally.
The EO problem minimizes 
\begin{equation}
    J = \frac{1}{2}\frac{T_{max}}{m_0}  \int^{t_1}_{t_0}  \Gamma_e^2 \  \textrm{d}t,  
    \label{3}
\end{equation}

subject to the dynamics in Eq. \eqref{d1}. The subscript $e$ is used to indicate quantities in the EO problem. The terminal constraints are the same as in Eq. \eqref{conditions2}.
The Hamiltonian of the EO problem, assuming $\lambda_m = 0$ at all times is
\begin{equation}
     \mathcal{H} =  \boldsymbol{\lambda_x}^T\left[ \boldsymbol{A} +   \boldsymbol{B} \frac{T_{max}}{m_0} \Gamma_e \widehat{\boldsymbol{\alpha}}_e\right] +  \frac{1}{2}  \frac{T_{max}}{m_0}\Gamma_e^2.
\end{equation}

The problem dynamics and the costate differential equations are again given by 
 {$\boldsymbol{\dot{x}} = \frac{\partial{\mathcal{H}}}{\partial{ \boldsymbol{\lambda}}} $ and $\boldsymbol{\dot{\lambda}}_{\boldsymbol{x}} = -\frac{\partial{\mathcal{H}}}{\partial{ \boldsymbol{x}}}$.} 
The thrust magnitude and direction that minimizes the Hamiltonian are $   \widehat{\boldsymbol{\alpha}}_e^* = -\frac{ \boldsymbol{B}^T \boldsymbol{\lambda_x}}{\|\boldsymbol{B^T \lambda_x}\| } \ \text{and} \  \Gamma_e^* = \|\boldsymbol{B^T \lambda_x}\|.$
Then, the optimal control problem   {yeilds a TPBVP similar in form to Eq. \eqref{Phi1}.}
This problem can be solved using a nonlinear solver under the boundary conditions given. The initial guess for solving the EO problem is obtained by linearizing it and solving for the initial costates analytically, as shown in \cite{Altitude-dependent2010c}.  \par 
\subsection{Calculating \texorpdfstring{$\Gamma_{TR}$}{GTR}}\label{Gammatr}
It can be noted that the FO switching function (Eq. \eqref{FOSF}) evaluated on the EO first guess is $\rho = \Gamma_{TR} - \Gamma_e$. Hence on the first guess

\begin{equation}
\begin{aligned}
\text {If}  \ \Gamma_e > \Gamma_{TR} &\Rightarrow \Gamma_f^* = \frac{m_0}{m(t)}  \\
\text {If} \  \Gamma_e < \Gamma_{TR} &\Rightarrow \Gamma_f^* = 0.  \\
     \end{aligned}
\end{equation}
Thus, we aim to compute the value of $\Gamma_{TR}$ such that the first guess control profile of the FO problem produces the same $\Delta v$ of the EO solution. The steps to calculate $\Gamma_{TR}$ are given in Algorithm \ref{trcalculation}.

\begin{algorithm}[hbt!]
\caption{ {Compute $\Gamma_{TR}$}}\label{trcalculation}
\begin{algorithmic}

\State Obtain the EO solution and its profile of $\Gamma_e(t)$ and calculate $\Delta v = \frac{T_{max}}{m_0}\int^{t_1}_{t_0} \Gamma_e(t) \ dt$.
\State Let $\Gamma_{TR_L}  = 0$ and $\Gamma_{TR_U}  = \max 
    \left(\Gamma_e(t) \right)$, set $\Delta v_{\Gamma_{TR}}  = 0$. 
\While{$| \Delta v - \Delta v_{\Gamma_{TR}}|  > tol $}
    \State Calculate $\Gamma_{TR} = \frac{\Gamma_{TR_L}+\Gamma_{TR_U}}{2}$.
 
     \If{$\Gamma_e(t) > \Gamma_{TR}$}    \Comment{Construct the bang bang thrust profile $\Gamma_{f}(t)$ from the continuous $\Gamma_e(t)$}
    \State $\Gamma_{f}(t)  = \frac{m_0}{m(t)} $   \Comment{ $m(t)$ is provided by Eq \eqref{m1}.}
  \Else
    \State $\Gamma_{f}(t) = 0$
  \EndIf
  
    \State Calculate $\Delta v_{\Gamma_{TR}} = \frac{T_{max}}{m_0}\int^{t_1}_{t_0}  \Gamma_{f}(t) \ dt$. 
      \If{$ \Delta v_{\Gamma_{TR}} > \Delta v$} \Comment{Apply the update rule}
    \State $\Gamma_{TR_L}  = \Gamma_{TR}$
  \Else
    \State $\Gamma_{TR_U}  = \Gamma_{TR}$
  \EndIf
    \EndWhile
    \State Return $\Gamma_{TR}$.
\end{algorithmic}
\end{algorithm}

In Fig. \ref{1c} of the following section, profiles of $\rho$ and $\Gamma$ are plotted against time to illustrate the effectiveness of determining $\Gamma_{TR}$ through this method. From this figure, it can be seen that setting $\Gamma_{TR}$ as described ensures that the profile of thrust provides an improved initial guess for the FO thrust profile.

\subsection{Overall Solution Methodology}

Algorithm \ref{method} summarises the solution methodology for solving FO problems.

\begin{algorithm}[hbt!]
\caption{Solve FO problem}\label{method}
\begin{algorithmic}
\State  {Solve the EO problem to obtain $\boldsymbol{\lambda}_{e}$ and calculate $\Gamma_{TR}$ as in Section \ref{Gammatr}.}
\State Let $i = 0, k = 0$ and guess solution $\boldsymbol{\lambda}_{g} = \boldsymbol{\lambda}_{e}$. 
\While{$k \leq k_{max}$}

    \State Solve the SFO problem with $\boldsymbol{\lambda}_{g}$ as the initial guess, to obtain $\boldsymbol{\lambda}(t_0)$.
    \State $\boldsymbol{\lambda}_{g} = \boldsymbol{\lambda}(t_0)$ , $k = k + \Delta k $ , $i = i+1$ 
    \EndWhile
    \State Solve the FO problem with $\boldsymbol{\lambda}_{g}$ as an initial guess.
\end{algorithmic}
\end{algorithm}

\section{Time-Optimal Problem}
The TO problem entails minimizing

\begin{equation}
    J = \beta_t \int^{t_1}_{t_0} t \  \textrm{d} t 
    \label{toptJ}
\end{equation}
subject to the dynamics given in Eq. \eqref{d1}. $t_0$ and $t_1$ denote the start and end times of the mission respectively and $\beta_t$ is a constant scaling factor that shall be determined later. The subscript $t$ is used to indicate quantities for in the TO problem. The terminal constraints are  
\begin{equation}\label{conditionst}
\boldsymbol{x}(t_0) = \boldsymbol{x_0}, \  \boldsymbol{x}(t_1) = \boldsymbol{x_1}, \ H(t_1) - \dot{L}_{t}\lambda_L(t_1) = 0,
\end{equation} 
in which $\dot{L}_{t}$ is the time derivative of the target true longitude. 
The TO problem can also be transformed into a TPBVP using Pontryagin’s maximum principle. First, the Hamiltonian is built as
\begin{equation}\label{Ht}
     \mathcal{H} =  \boldsymbol{\lambda_x}^T\left[ \boldsymbol{A} + \boldsymbol{B} \frac{T_{\max }}{m_0}\Gamma_t \widehat{\boldsymbol{\alpha}}_t\right] +  \beta_t,
\end{equation}
The $ \widehat{\boldsymbol{\alpha}}_t$ and $\Gamma_t$ that minimize the Hamiltonian are  {$
    \widehat{\boldsymbol{\alpha}}_t^* = -\frac{ \boldsymbol{B^T \lambda_x}}{\|\boldsymbol{B^T \lambda_x}\|} \ \text{and} \  \Gamma_t^* = \frac{m_0}{m(t)}.$} The initial and final states 
    must satisfy conditions in Eq. \eqref{conditionst}. Thus, the problem now consists of finding $\boldsymbol{\lambda_{x}}(t_0)$ that satisfy 
\begin{equation}\label{Phit}
    \boldsymbol{\Phi}(\boldsymbol{\lambda_{x}}(t_0)) = [ \boldsymbol{x}(t_1) - \boldsymbol{x}_1, \mathcal{H} (t_1) - \dot{L}_t \lambda_L(t_1)]^T = \boldsymbol{0}.
\end{equation}
This shooting function can then be solved using a nonlinear solver.
\par 
If the spacecraft mass were treated as a state variable, the expression of the Hamiltonian would have a term with the mass costate ($\lambda_m$). Hence $\mathcal{H} (t_1)$ in Eq. \eqref{Phit} would be dependent on $\lambda_m(t_1)$. However, as there are no constraints on the final mass, $\lambda_m(t_1)$  must be zero, and, consequently, the dynamics of $\lambda_m$ do not impact the solution. For this reason, the dynamics of $\lambda_m$ can be ignored. 

\subsection{Initial Guess for the Optimal Time of Flight} \label{tofguess}
A non-intuitive initial guess is required for the time of flight of the TO solution. The EO problem described in section \ref{eoptsection} is used for this purpose. The procedure to obtain a guess for the optimal time of the transfer is given in Algorithm \ref{a3}.

\begin{algorithm}[hbt!]
\caption{ {Generate guess for the time of flight}}\label{a3}
\begin{algorithmic}

\State Let $t_{L} = 0$,  $t_{U} =  \max(TOF)$, $\Delta v_{t} =  0$  $\Delta v_{e} =  1$. 
\While{$|\Delta v_{t} - \Delta v_{e}|  > tol $}
    \State Calculate $t =  \frac{t_{L}+t_{U}}{2}$.
    \State Obtain the EO solution and $\Gamma_e$ for time of flight $t$. 
    \State Calculate $\Delta v_{e}=  \frac{T_{max}}{m_0} \int_{0}^{t}  \Gamma_e 
        \ dt $  and $\Delta v_{t}= \frac{T_{max}}{m_0} \int_{0}^{t}  \frac{m_0}{m(t)} dt $.

      \If{$ \Delta v_{t} > \Delta v_{e}$} \Comment{Apply the update rule}
    \State ${t_U} = t$
  \Else
    \State $ {t_L} = t$
  \EndIf
    \EndWhile
    \State Return $t$.
\end{algorithmic}
\end{algorithm}
With this procedure we determine the time of flight ($t$) as the one for which the $\Delta v$ of the EO problem is equal to that obtained thrusting all time (as in the TO solution) for that time.
 \subsection{Calculating \texorpdfstring{$\beta_t$}{bt} }\label{calcbetat}
 Once we have obtained the EO solution for the new time of flight $t$, we calculate $\mathcal{H}(t_1) - \dot{L}_t\lambda_L(t_1) $ for the first iteration of the solver. Then, we calculate $\beta_t$ by setting  $\mathcal{H}(t_1) - \dot{L}_t\lambda_L(t_1) =0 $. This ensures that the constraint is satisfied on the initial guess. The complete procedure for calculating $\beta_t$ is given in Algorithm \ref{a4}.
 \begin{algorithm}[hbt!]
\caption{ {Compute \texorpdfstring{$\beta_t$}{bt}}}
\label{a4}
\begin{algorithmic}
\State  Solve EO problem to obtain $\boldsymbol{\lambda}_{e}(t_0)$.
\State Propagate the TO dynamics with $\boldsymbol{\lambda}_{e}(t_0)$ for the same time period as the EO. 
\State  Calculate $\widehat{\boldsymbol{\alpha}}_{t_1}$ and $\Gamma_{t_1}$ as
     \begin{equation}
         \widehat{\boldsymbol{\alpha}}_{t_1} = -\frac{ {\boldsymbol{B}^T \boldsymbol{\lambda}_{e}(t_1)}}{\|{\boldsymbol{B}^T \boldsymbol{\lambda}_{e}(t_1)}\|} \ \text{and} \ \Gamma_{t_1}  = \frac{m_0}{m(t_1)}.
     \end{equation}
\State Calculate $\mathcal{H}(t_1) - \dot{L}_t \lambda_L(t_1) $ using Eq. \eqref{Ht} as 
      \begin{equation}
     \mathcal{H}(t_1) - \dot{L}_t \lambda_L(t_1) = \boldsymbol{\lambda_{e}}(t_1)^T\left[ \boldsymbol{A}+ \boldsymbol{B} \frac{T_{\max}}{m_0} \Gamma_{tf} \widehat{\boldsymbol{\alpha}}_{t_1}\right] + \beta_t -\dot{L}_t\lambda_L(t_1). 
 \end{equation}
 \State Set $\mathcal{H}(t_1) - \dot{L}_t \lambda_L(t_1)  = 0$ to determine $\beta_t$.
  \begin{equation}\label{betat}
    \beta_t = \dot{L}_t\lambda_L(t_1)  - \boldsymbol{\lambda}_{e}(t_1)^T\left[ \boldsymbol{A} + \boldsymbol{B} \frac{T_{\max}}{m_0} \Gamma_{t_1} \widehat{\boldsymbol{\alpha}}_{t_1}\right]
 \end{equation}

\end{algorithmic}
\end{algorithm}

The impact of the value of $\beta_t$ on the number of iterations to solve a TO problem is shown in Fig. \ref{betafig}. The result for $\beta_t$ set by Eq. \eqref{betat} is represented  {by a star marker}. This figure shows that setting $\beta_t$ such that $\mathcal{H}(t_1) - \dot{L}_t \lambda_L(t_1)  = 0$ at the first iteration helps reducing the number of iterations required to solve the TO problem.
 
 \subsection{Overall Solution Methodology}

Algorithm \ref{method2} summarizes the full solution methodology discussed for solving time-optimal problems. 
\begin{algorithm}[hbt!]
\caption{Solve TO problem}\label{method2}
\begin{algorithmic}
\State Generate an initial guess for the time of flight as discussed in Section \ref{tofguess}. 
\State Solve the EO problem for the time of flight guess to obtain $\boldsymbol{\lambda}_{e}$.
\State Determine $\beta_t$ as discussed in Section \ref{calcbetat}.
\State Solve the TO problem using $\boldsymbol{\lambda}_{e}$ as an initial guess.
\end{algorithmic}
\end{algorithm}
\section{Case Studies}
 In this section, fuel and time-optimal trajectories are designed for several missions studied previously in the literature. All simulations are executed on a desktop computer with an Intel Core i5 CPU of 2.4Hz and 8.00 GB of RAM. The program is written in C++ and compiled with Microsoft Visual Studio Code. Simple shooting and the GNU Scientific Library's multidimensional root finding solver \cite{galassi2002gnu} are used to generate the following solutions. \par
Canonical units are used in the numerical simulations, such that the length unit is set to $L_u = 1$ AU (149,597,870.66 km) for heliocentric cases, and $L_u = 1$ $R_E$ (6378.1363 km) for geocentric cases. The time unit is set to  $T_u = 1$ year ($86400 \times 365.25$ s) for heliocentric cases and to $T_u = 1$ day ($86400$ s) for geocentric cases. The mass unit is set to the initial mass of spacecraft ($M_u = m_0$) in both heliocentric and geocentric cases.\par   {Table \ref{propsys} gives the spacecraft parameters used for each mission. The propulsion systems used in the interplanetary cases are compatible with gridded ion thrusters such as NASA's Evolutionary Xenon Thruster \cite{nasa_2021}. These thrusters have an efficiency of 0.7, corresponding to input powers of 12.61 kW and 6.73 kW, respectively, for Earth to Tempel 1 and Earth to Dionysus test cases. It is worth mentioning that these values were selected to ease the comparison of our approach with others available in the literature rather than from system/mission considerations.
The propulsion system used for the geocentric case is compatible with green chemical propulsion systems such as the B20 Thruster from Dawn Aerospace \cite{satsearch}.}
\begin{table}[hbt!]
\centering 
\caption{ {Propulsion system and spacecraft parameters for each mission}}\label{propsys}
\begin{tabular}{lccc}\toprule 
Mission    & $I_{sp}$ (s) & $T_{max}$ (N) & $m_0$ (kg)  \\\hline 
Earth to Tempel 1  &  3000  & 0.6    & 1000    \\
Earth to Dionysus   & 3000  & 0.32   & 4000      \\
Space debris & 300   & 1.0      & 100         \\\bottomrule       
\end{tabular}
\end{table}
\subsection{Earth to Tempel 1}
Given the launch time, the initial ($\boldsymbol{x}_0$) and target ($\boldsymbol{x}_1$) coordinates of this heliocentric transfer problem are computed using the Jet Propulsion Laboratory Horizons system\footnote{Data available online at \url{http://ssd.jpl.nasa.gov/?horizons} [retrieved 1 August 2021]} and  {are $\boldsymbol{x}_0 = [1.000064, -0.003764,  0.015791, \SI{-1.211e-5}{}, \SI{-4.514e-6}{}, 5.51356]^T$ and $\boldsymbol{x}_1 = [2.328616, -0.191235, -0.472341, 0.033222, 0.085426, 4.96395]^T$ in MEE coordinates}.  {The propulsion parameters and  boundary states are same as those given in \cite{Altitude-dependent2010c}. } 
\subsubsection{Fuel-Optimal Solution}
\begin{figure}[!hbt]
    \centering
    \subfloat[\centering Optimal Trajectory  \label{1a}]{{\includegraphics[width=8.0cm]{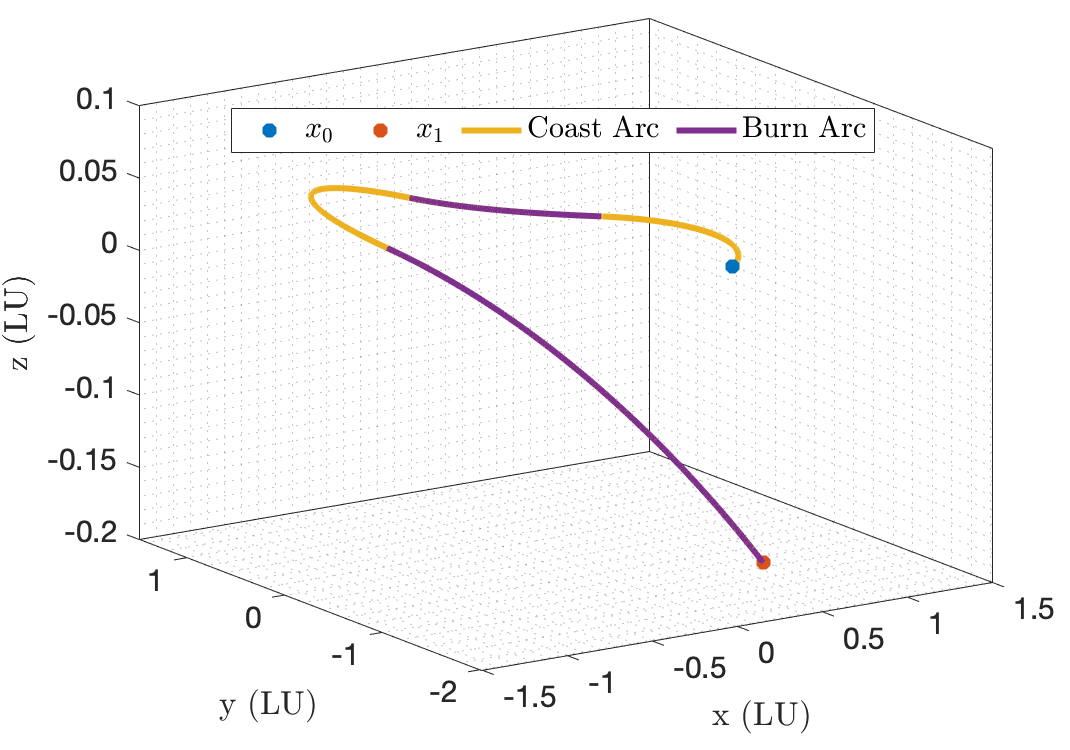} }}%
    \qquad
    \subfloat[\centering Control and spacecraft mass profiles \label{1b}]{{\includegraphics[width=7.5cm]{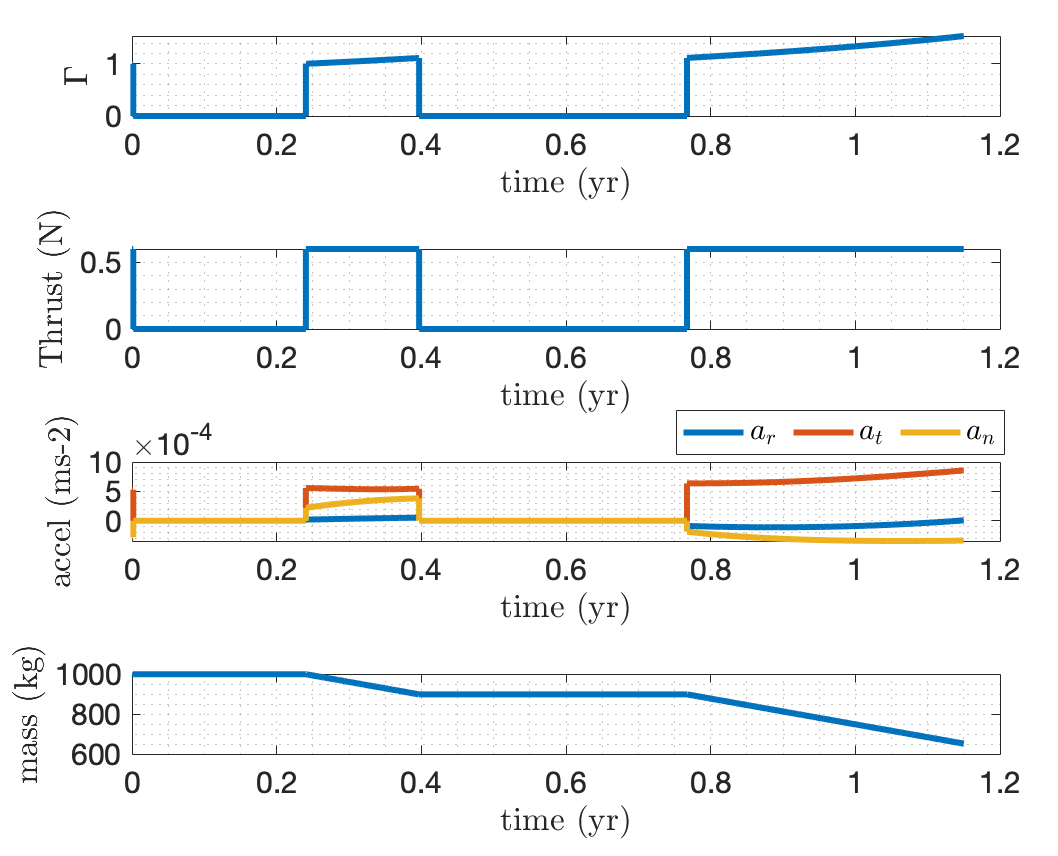} }}%
    \caption{Fuel-optimal Earth to Tempel 1 transfer}%
   
\end{figure}

 {Tables \ref{costatesinit1} and \ref{tsol} show the FO results obtained for a time of flight of 420 days.} Five continuation steps were used, and $k_{max}$ was taken to be 0.99.  Figure \ref{1a} shows the optimal trajectory, which contains two burn arcs and two coast arcs. Fig. \ref{1b} shows the optimal profile of thrust, $\Gamma$, acceleration, and mass consumption of the fuel optimal solution. \par 
This problem is also solved using a switched method in \cite{Altitude-dependent2010c}. Five continuation steps are taken to derive the same solution with a reported computational time of 0.234 s. Another solution to this problem is derived in \cite{Taheri2018} using a hyperbolic tangent smoothing method. It is reported to take five continuation steps and a computational time of 0.263 s. In \cite{Taheri2018}, the initial costates are randomly generated within a specified range. Hence, the convergence rate of this method is determined by sampling many random guesses. Instead, our method and the switched method \cite{Altitude-dependent2010c} calculate the initial costates from a linearized EO problem, so they both eliminate the need for random guess generation. For this simple problem, our method, the switched method, and the hyperbolic tangent smoothing method all have a 100\% convergence rate \cite{Altitude-dependent2010c}. \par
\begin{figure}[!hbt]
    \centering
        \subfloat[\centering Profiles of $\rho$ and $\Gamma$ \label{1c} ]{{\includegraphics[width=6.5cm]{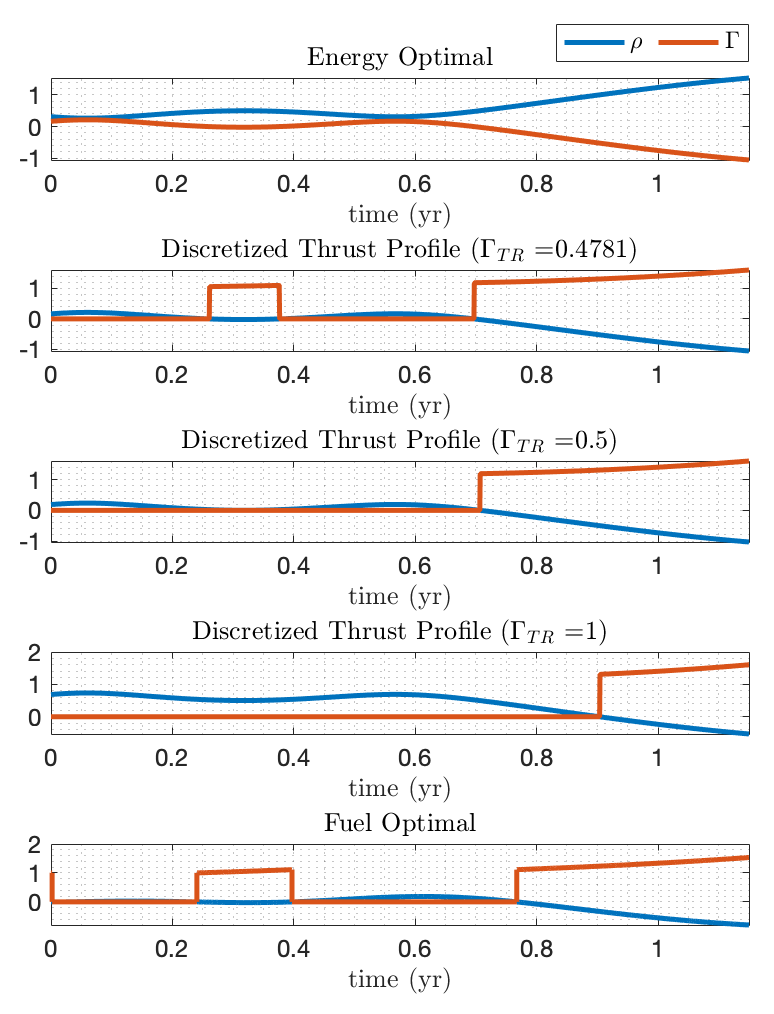} }}%
        \qquad 
    \subfloat[\centering Profile of mass at each step of the solution process (top) and the profile of $\Gamma$ at each step of the solution process (bottom) \label{1d}]{{\includegraphics[width=7.3cm]{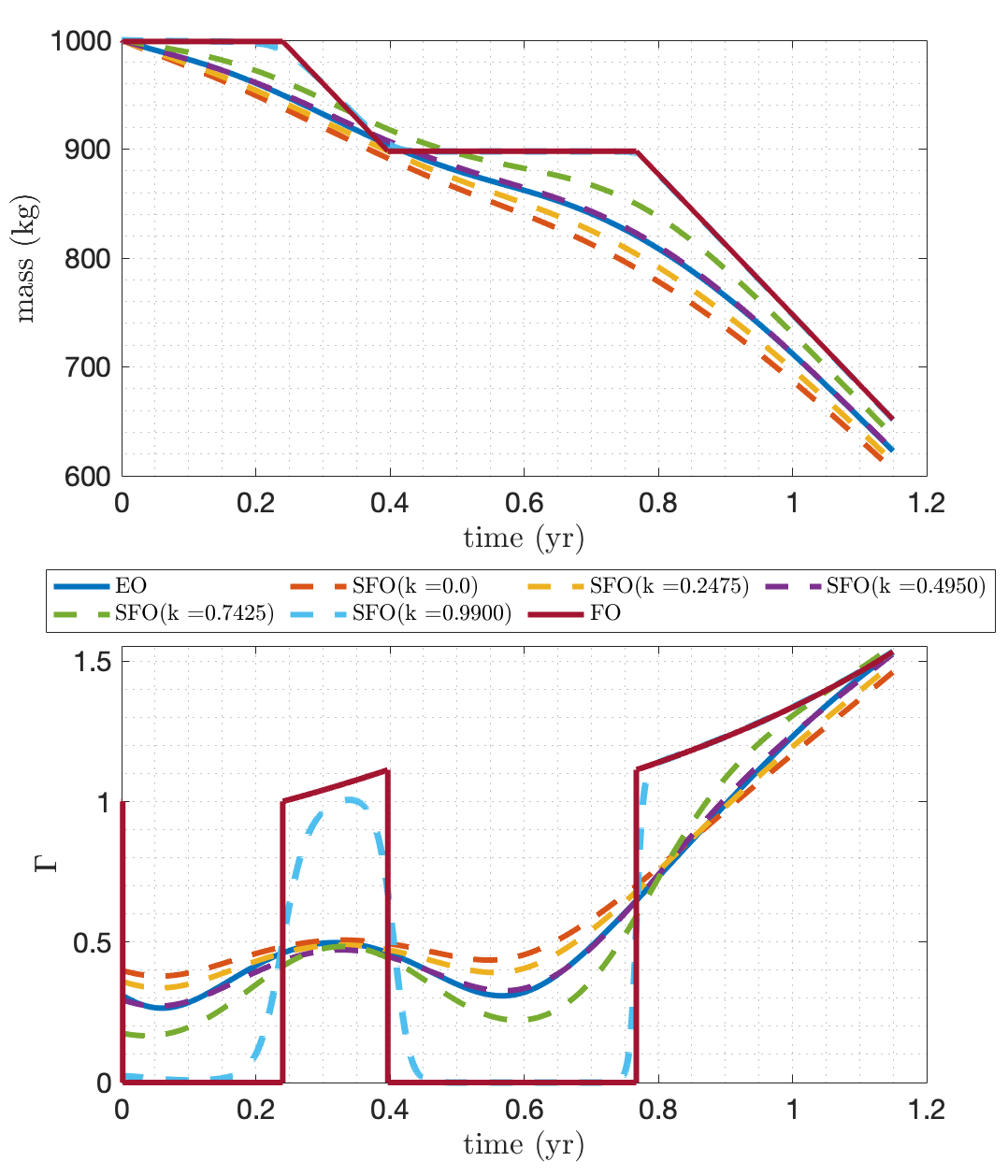} }}%
     \qquad
         \subfloat[\centering Angle between the thrust direction at each step of the solution process and the optimal thrust direction. \label{1e}]{{\includegraphics[width=7.0cm]{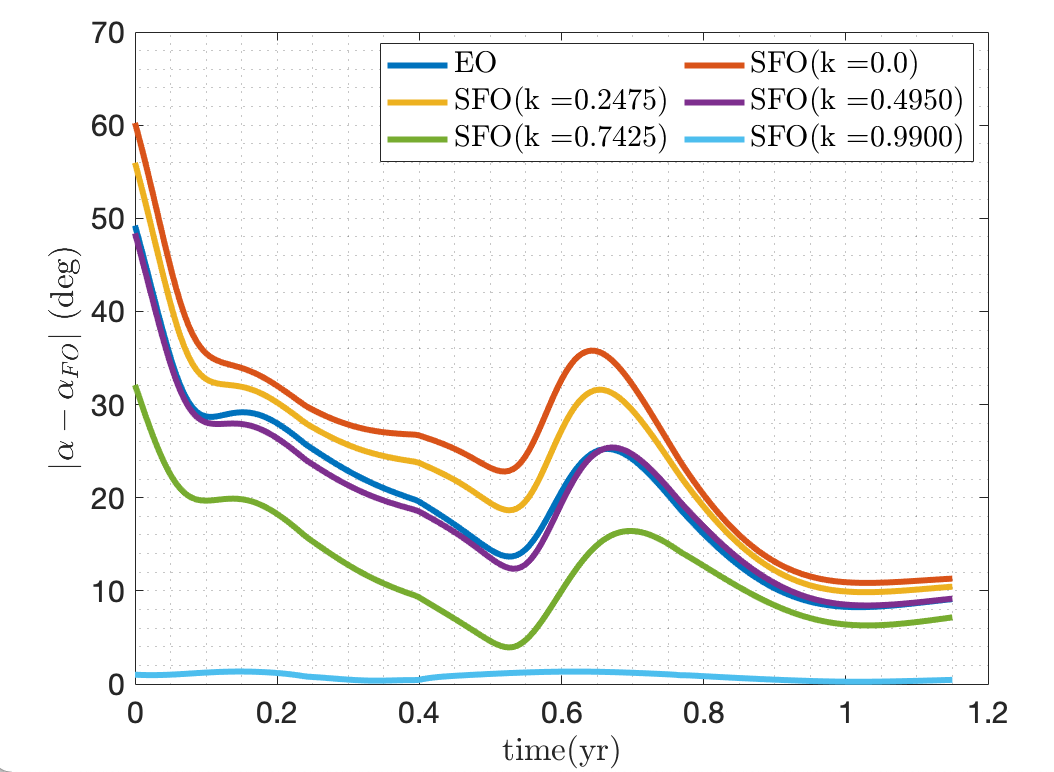} } }
        \subfloat[\centering Profiles of the switching function ($\rho$) \label{1f} ]{{\includegraphics[width=6.5cm]{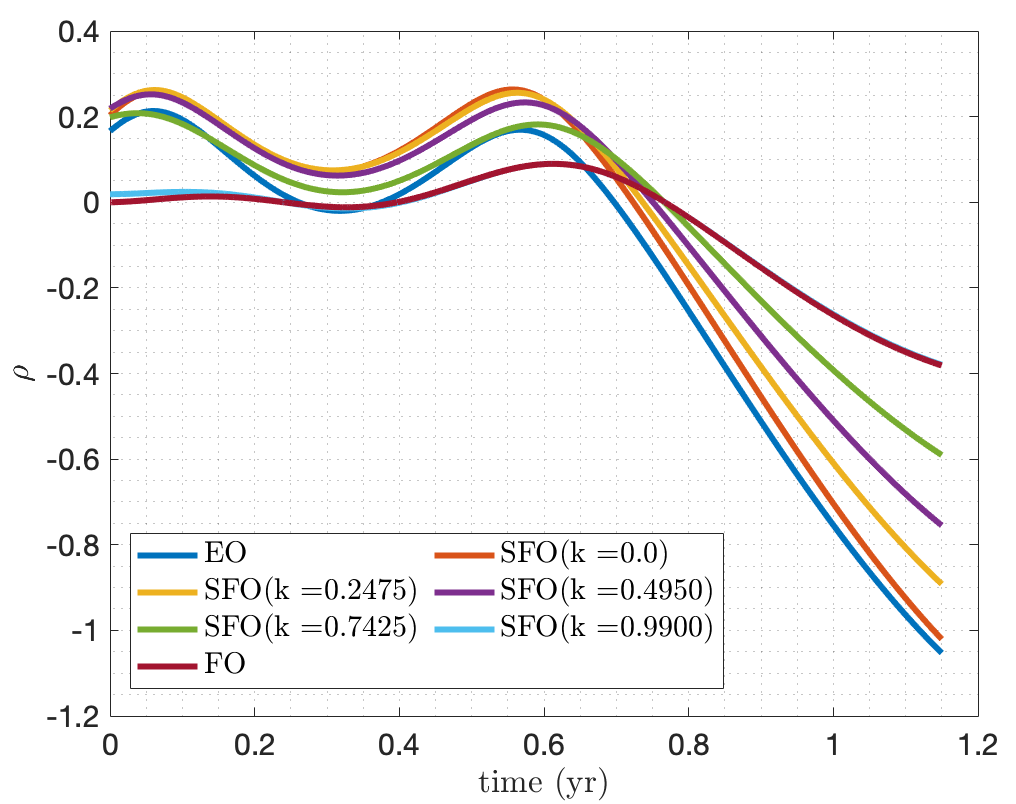} }}%
             \qquad
   \caption{Solution generation process of Earth to Tempel 1 trajectory, FO case}%
\end{figure}

 
Figure \ref{1c} shows the change in $\rho$ and $\Gamma$ for different $\Gamma_{TR}$ values.  $\Gamma_{TR} =  0.4781$ is the case where the value of $\Gamma_{TR}$ is obtained from the thresholding method described. It can be seen that the profile of $\rho$ generated with the $\Gamma_{TR}$ value from the thresholding method matches very closely to the $\rho$ profile of the final FO solution. Note that for $\Gamma_{TR} =1$ (i.e., without using any scaling), the initial switching function profile is very different from the one corresponding to the optimal solution, indicating that the thresholding does indeed provide a better initial guess. Similarly, it can be seen that the $\Gamma_{TR}$ value generated from the thresholding method produces a $\Gamma$ that matches very closely with the final FO solution obtained. This feature enables faster convergence when the FO is solved.  {When collocation methods are used to solve the problem, the thrust profile matching through $\Gamma_{TR}$ provides an immense advantage and aids convergence. However, in simple cases, this is even beneficial when single shooting methods are used. } \par
Figure \ref{1d} shows the profiles of $\Gamma$ and mass. The total propellant consumption of the mission agrees with the values reported in \cite{Altitude-dependent2010c, Taheri2018}. Figure \ref{1e} shows how the optimal thrust direction compares to the optimal thrust direction of the FO problem throughout the solution process. The final SFO thrust direction matches closely with that of the FO solution. In Fig. \ref{1f} the variation of the switching function going from EO to FO is reported, showing that the switching function of the EO provides a good initial guess for that of FO.

\begin{table}[!hbt]
\centering 
\caption{Initial costates obtained at each step of the solution for the fuel optimal Earth to Tempel 1 trajectory}
\begin{tabular}{llccccccc}
\toprule
Method & Dynamic & Fuel mass & $\lambda_p$ & $\lambda_f$& $\lambda_g$ & $\lambda_h$  & $\lambda_k$& $\lambda_L$\\
~  & Model              &kg & ${L_u/T_u^3}$ & ${L_u^2/T_u^3}$& ${L_u^2/T_u^3}$ & ${L_u^2/T_u^3}$  & ${L_u^2/T_u^3}$& ${L_u^2/T_u^3}$ \\ \hline
EO              & Kep.    & 377.2121  & 0.5554        & -1.5382         & -0.3929         & -1.2909         & -5.0413         & -0.4974         \\
SFO (k = 0.00)      & Kep.    & 394.6693  & 0.8148        & -1.6150         & -0.3390         & -1.3274         & -4.3093         & -0.5047         \\
SFO (k = 0.2475) & Kep.    & 387.0673  & 0.6535        & -1.4499         & -0.3352         & -1.2074         & -4.1986         & -0.4547         \\
SFO (k = 0.4950)  & Kep.    & 376.6296  & 0.4238        & -1.2649         & -0.3468         & -1.0543         & -4.1247         & -0.3969         \\
SFO (k = 0.7425) & Kep.    & 363.3607  & 0.0105        & -1.0190         & -0.3592         & -0.8020         & -3.9924         & -0.3128         \\
SFO (k = 0.99) & Kep.    & 348.5101  & -0.8875       & -0.5664         & -0.0808         & 0.1927          & -3.6568         & -0.1386         \\
FO              & Kep.    & 348.2554  & -0.9249       & -0.5600         & -0.0446         & 0.2963          & -3.6778         & -0.1315         \\ \bottomrule
\end{tabular}
\label{costatesinit1}
\end{table}

\begin{figure}[hbt!]
    \centering
    \subfloat[\centering The effectivity of calculating $\Gamma_{TR}$ in FO case  \label{GammaTRiter}]{{\includegraphics[width=7.6cm]{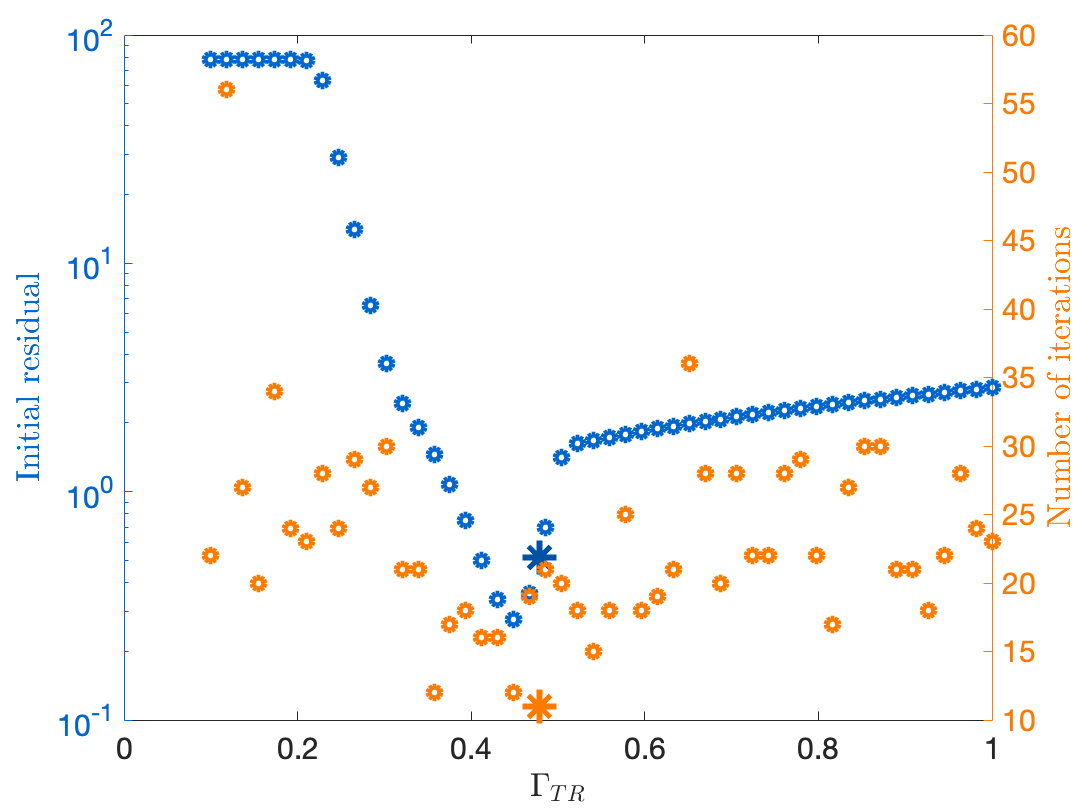} }}%
    \qquad
    \subfloat[\centering The effectivity of calculating $\beta_t$ in TO case \label{betafig}]{{\includegraphics[width=7.6cm]{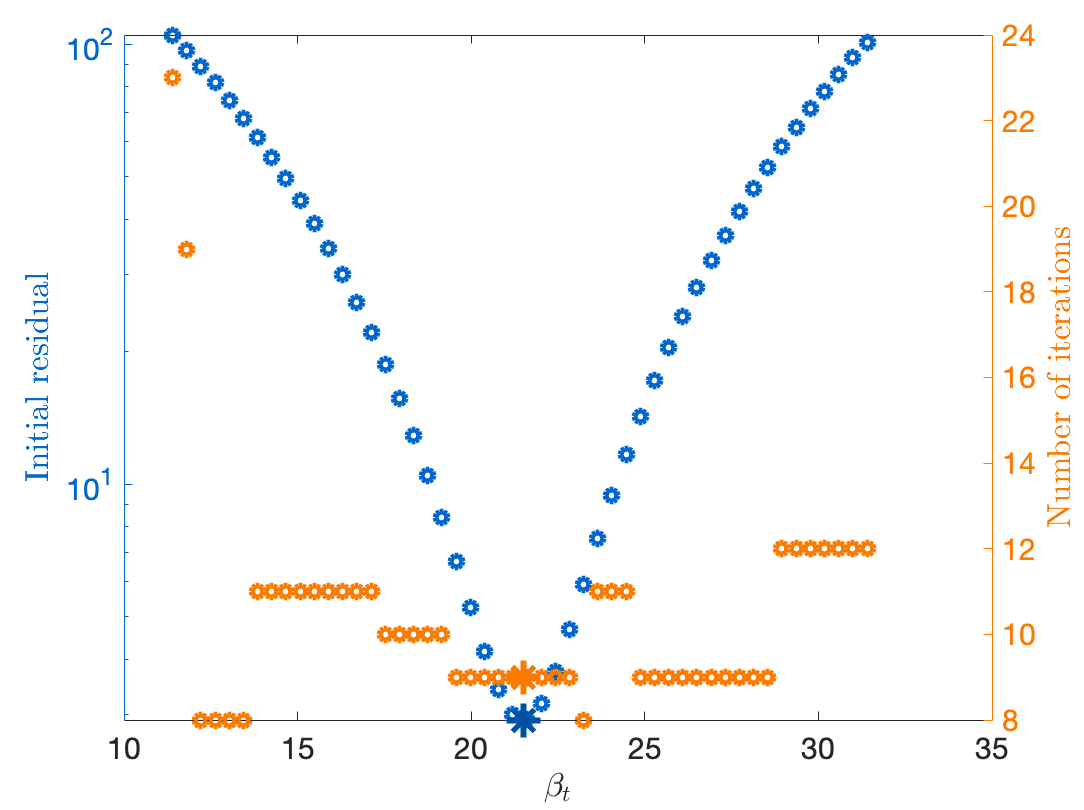} }}%
    \caption{ {Effectivity of $\Gamma_{TR}$ and $\beta_t$ calculations for Earth to Tempel 1 case studies. (Note that the values for the calculated $\beta_t$ and $\Gamma_{TR}$ are indicated by star markers) }}%

\end{figure}
Note that while the smoothing is required for solving this problem through the multidimensional root finding solver, it can be skipped if more nuanced algorithms such as Matlab's trust-region or trust-region-dogleg algorithms are used.  \par 
 {Figure \ref{GammaTRiter} shows how the number of iterations required to solve the FO (right y-axis) and the initial residual at the start of the continuation (left y-axis) vary with $\Gamma_{TR}$. The $\Gamma_{TR}$ set by the method discussed is shown by a star marker. This figure illustrates that by setting $\Gamma_{TR}$ this way, we can significantly reduce the number of iterations required by reducing the initial residual. While the initial residual obtained using the calculated $\Gamma_{TR}$ is not the minimum residual possible, it is still sufficiently low to result in the smallest number of iterations reported for convergence.}

\subsubsection{Time-Optimal Solution}
\begin{figure}[hbt!]
    \centering
    \subfloat[\centering Optimal Trajectory  \label{dt1a}]{{\includegraphics[width=7.5cm]{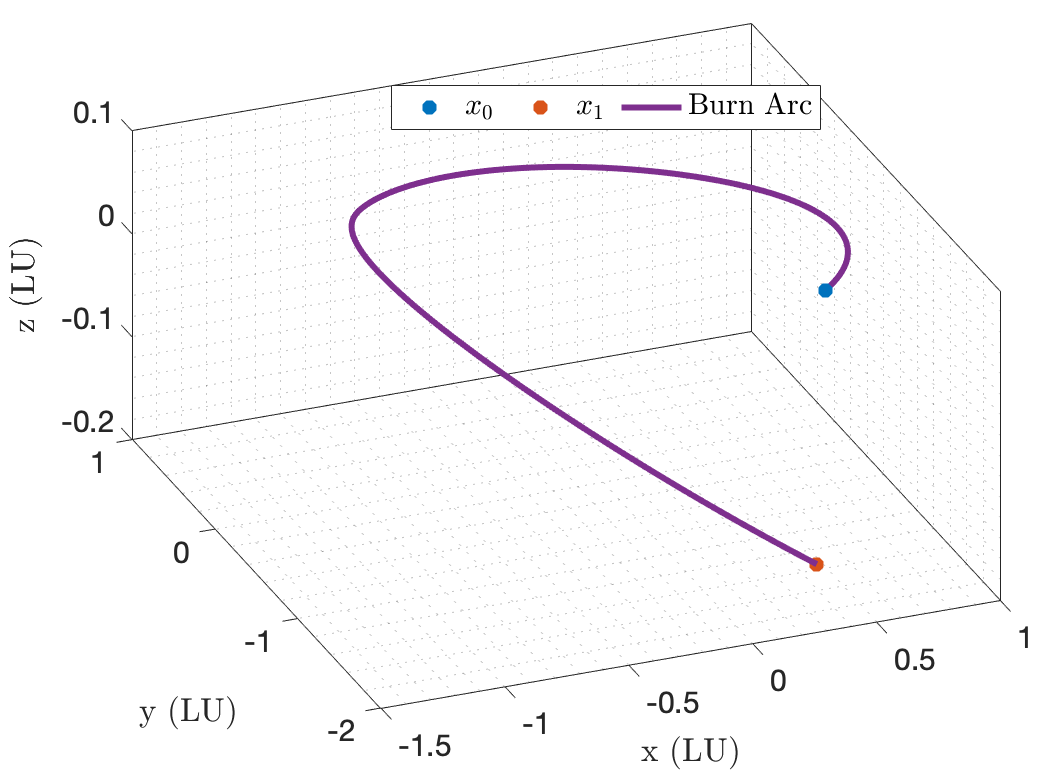} }}%
    \qquad
    \subfloat[\centering Control and spacecraft mass profiles \label{dt1b}]{{\includegraphics[width=7.5cm]{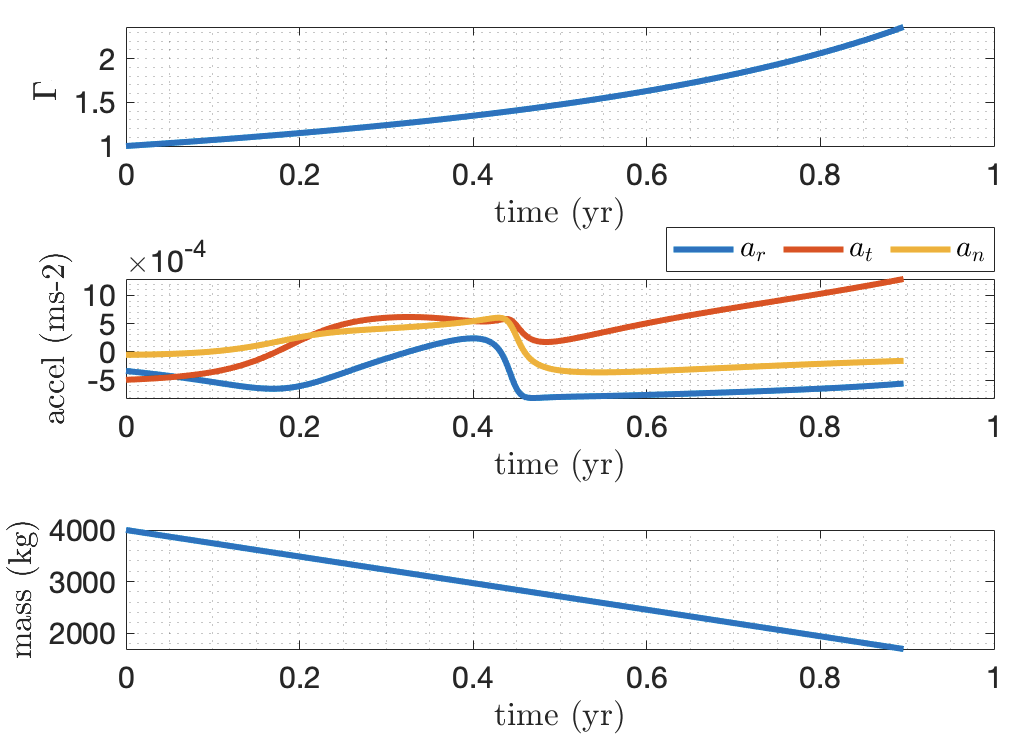} }}%
    \caption{Time-optimal Earth to Tempel 1 transfer}%
\end{figure}

 {Tables \ref{topt} and \ref{tsol}} show the TO results obtained. Fig. \ref{dt1a} shows the time-optimal trajectory while Fig. \ref{dt1b} show the profiles of acceleration, $\Gamma$, and spacecraft mass. \par 


\begin{table}[!hbt]
\centering 
\caption{Initial costates obtained at each step of the solution for the time-optimal Earth to Tempel 1 trajectory}
\begin{tabular}{llccccccc}
\toprule
Method                                                           & Dynamic& TOF                      & $\lambda_p$                       & $\lambda_f$                                                             & $\lambda_g$                         & $\lambda_h$                         & $\lambda_k$                         & $\lambda_L$                         \\
~  & Model              &days & ${L_u/T_u^3}$ & ${L_u^2/T_u^3}$& ${L_u^2/T_u^3}$ & ${L_u^2/T_u^3}$  & ${L_u^2/T_u^3}$& ${L_u^2/T_u^3}$ \\
\hline
EO & Kep. &  420.0000   & 0.5554  & -1.5382  & -0.3929  & -1.2909 &  -5.0413 &  -0.4974 \\
EO (for TOF guess) & Kep. &   307.7231 &  7.7238 &  -7.2184&   -5.0119 &  -1.7433 &  -6.5260 &  -1.8145    \\
TO  & Kep. &    327.1544&  3.0987 &  -2.7809  & -1.6998  & -0.3594  & -1.7373  & -0.7050 \\  \bottomrule
\end{tabular}
\label{topt}
\end{table}


 {Figure \ref{betafig} illustrates how the number of iterations required to solve the TO problem (right y-axis)  and the initial residual (left y-axis)  changes with $\beta_t$ for this test case. From this plot, it can be seen that setting $\beta_t$ as described in section \ref{calcbetat} results in a very low initial residual, which makes convergence easier. }


\subsection{Earth to Dionysus}

This test case is also adapted from \cite{Altitude-dependent2010c}.  The initial and target states are $\boldsymbol{x}_0 = [0.999316, -0.004023, 0.015873,  \SI{-1.623e-5}{}, \SI{1.667e-5}{}, 1.59491]^T $ and $\boldsymbol{x}_1 = [1.555261 , 0.152514, -0.519189, 0.016353, 0.117461, 2.36696]^T$ in MEE coordinates. The maximum thrust is lower than in the previous test case, while the spacecraft mass is higher. Hence the initial acceleration produced is approximately $1/8$ \textsuperscript{th} of the previous mission. 
\subsubsection{Fuel-Optimal Solution}
\begin{figure}[hbt!]
    \centering
    \subfloat[\centering Optimal Trajectory    \label{dionmethod1a}]{{\includegraphics[width=7.5cm]{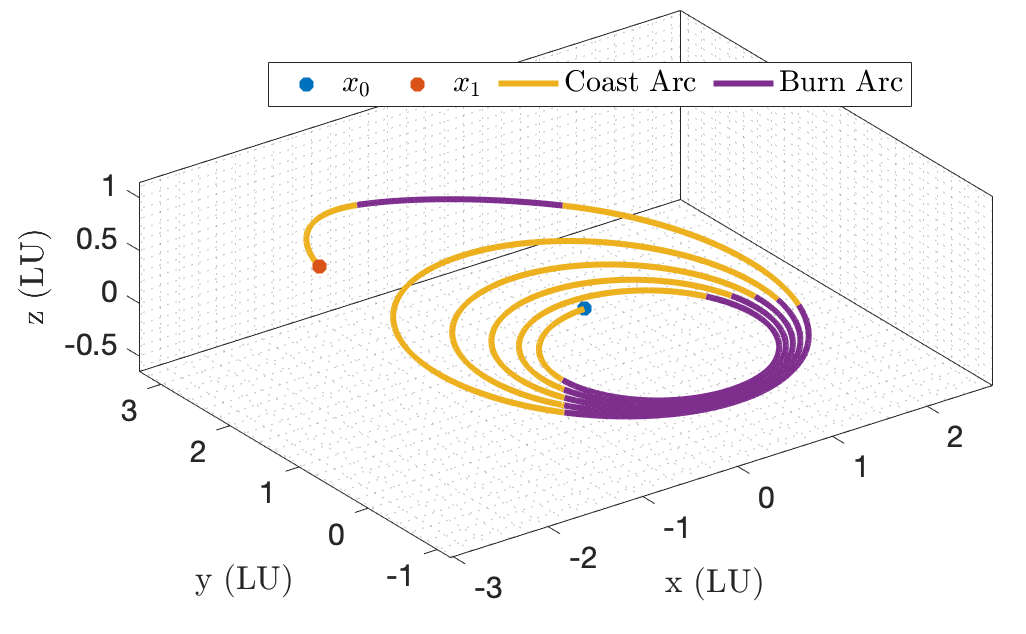} }}%
    \qquad
    \subfloat[\centering Control and spacecraft mass profiles    \label{dionmethod1b}]{{\includegraphics[width=7.5cm]{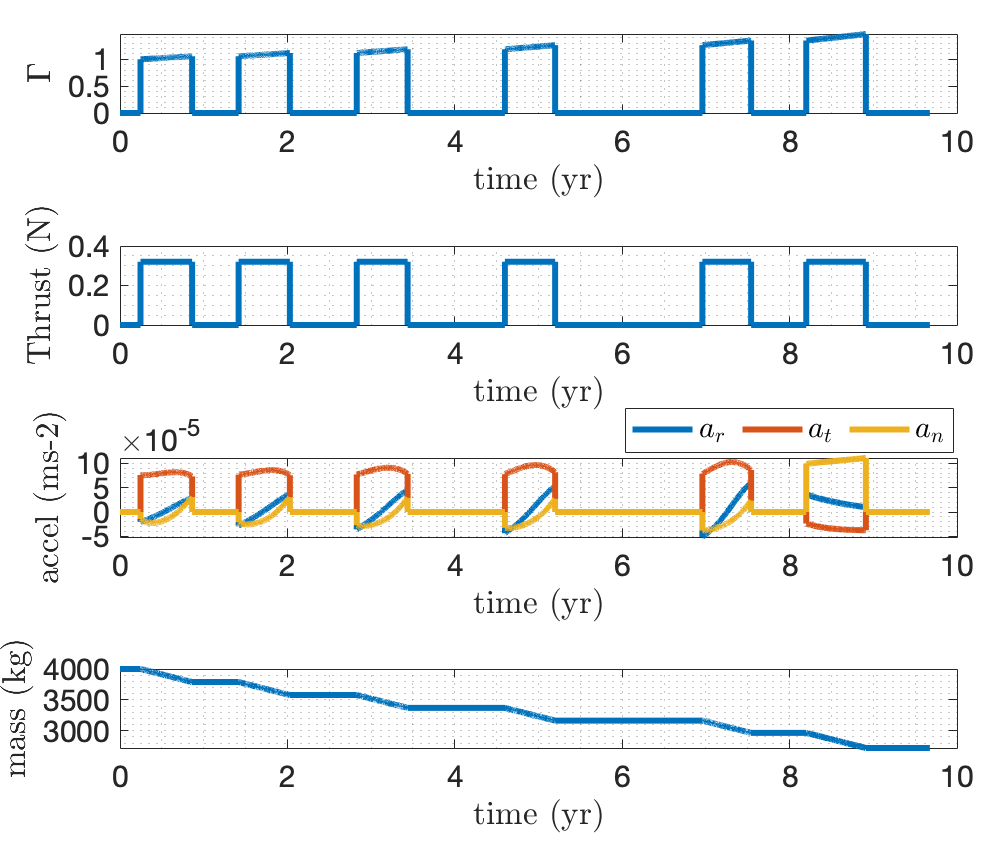} }}%
    \caption{Fuel-optimal Earth to Dionysus transfer}%
 
\end{figure}

 {Tables \ref{dioncost1} and \ref{tsol} show the FO results obtained, for a time of flight of 3534 days.} Again, five continuation steps are required to generate a solution, and $k_{max}$ was taken to be 0.99.   Figure \ref{dionmethod1a} shows the FO trajectory obtained, containing seven coast arcs and six burning arcs. Fig. \ref{dionmethod1b} illustrates the profiles of $\Gamma$ and mass of the spacecraft, as well as the trajectory taken. The proposed method has a convergence rate of 100\%.  \par
 {For the same problem, the switched method in \cite{Altitude-dependent2010c} takes five continuation steps to generate an optimal fuel consumption of 1279.93 kg within a reported computational time of 0.419 s. The hyperbolic tangent smoothing method  \cite{Taheri2018} is reported to take five continuation steps to generate a 1281.68 kg solution in 0.819 s. The switched method reports a 100\% convergence rate, while the smoothing method only has a reported convergence rate of 68.3\%.}\par 

\begin{table}[!hbt]
\centering 
\caption{Initial costates obtained at each step of the fuel optimal solution for the Earth to Dionysus trajectory}
\begin{tabular}{llccccccc}
\toprule
Method          & Dynamic & Fuel mass & $\lambda_p$   & $\lambda_f$     & $\lambda_g$     & $\lambda_h$     & $\lambda_k$     & $\lambda_L$     \\
~          & Model   & kg        & ${L_u/T_u^3}$ & ${L_u^2/T_u^3}$ & ${L_u^2/T_u^3}$ & ${L_u^2/T_u^3}$ & ${L_u^2/T_u^3}$ & ${L_u^2/T_u^3}$ \\ \hline
EO              & Kep.    & 1479.0246 & -1.7649       & -0.2215         & 1.0965          & -1.0684         & -2.3545         & -0.0096         \\
SFO (k =0.0)      & Kep.    & 1590.4344 & -1.7831       & -0.1666         & 0.8657          & -0.8845         & -1.5500         & -0.0027         \\
SFO (k =0.2475) & Kep.    & 1556.4834 & -1.8173       & -0.1715         & 0.8662          & -0.8952         & -1.6566         & -0.0006         \\
SFO (k =0.495)  & Kep.    & 1500.5678 & -1.8506       & -0.1833         & 0.8806          & -0.9193         & -1.8726         & 0.0017          \\
SFO (k =0.7425) & Kep.    & 1401.8552 & -1.8615       & -0.2100         & 0.9037          & -0.9478         & -2.3164         & 0.0040          \\
SFO (k =0.99) & Kep.    & 1281.1310 & -1.8944       & -0.2917         & 0.8610          & -1.1312         & -2.6800         & 0.0065          \\
FO              & Kep.    & 1280.7021 & -1.8944       & -0.2923         & 0.8609          & -1.1313         & -2.6767         & 0.0065          \\ \bottomrule
\end{tabular}
\label{dioncost1}
\end{table}

\subsubsection{Time-Optimal Solution}
\begin{figure}[hbt!]
    \centering
    \subfloat[\centering Optimal Trajectory   \label{dta}]{{\includegraphics[width=8.0cm]{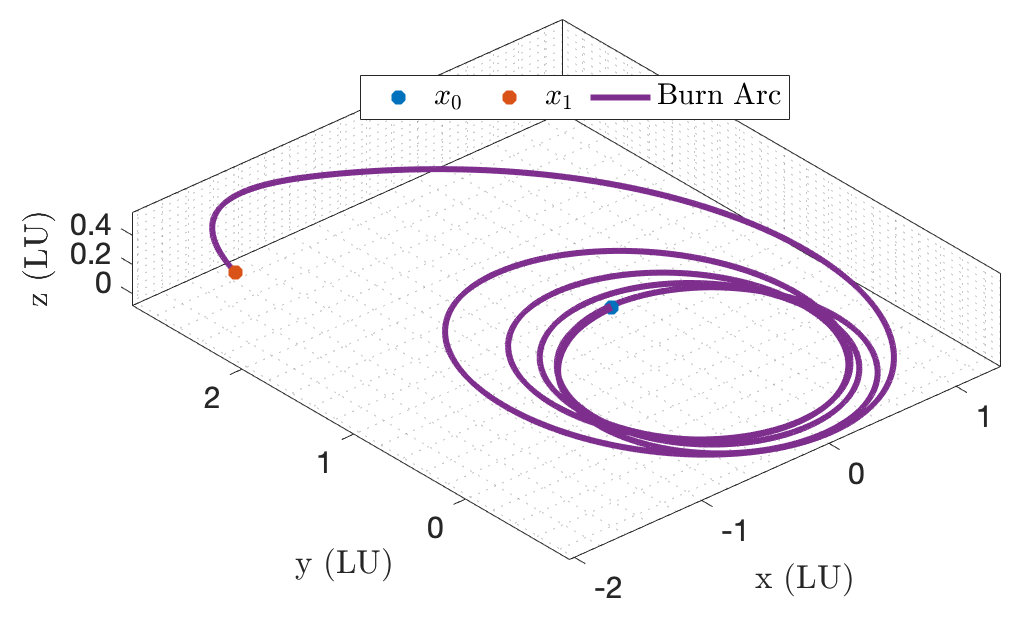} }}%
    \qquad
    \subfloat[\centering Control and spacecraft mass profiles \label{dtb}]{{\includegraphics[width=7.5cm]{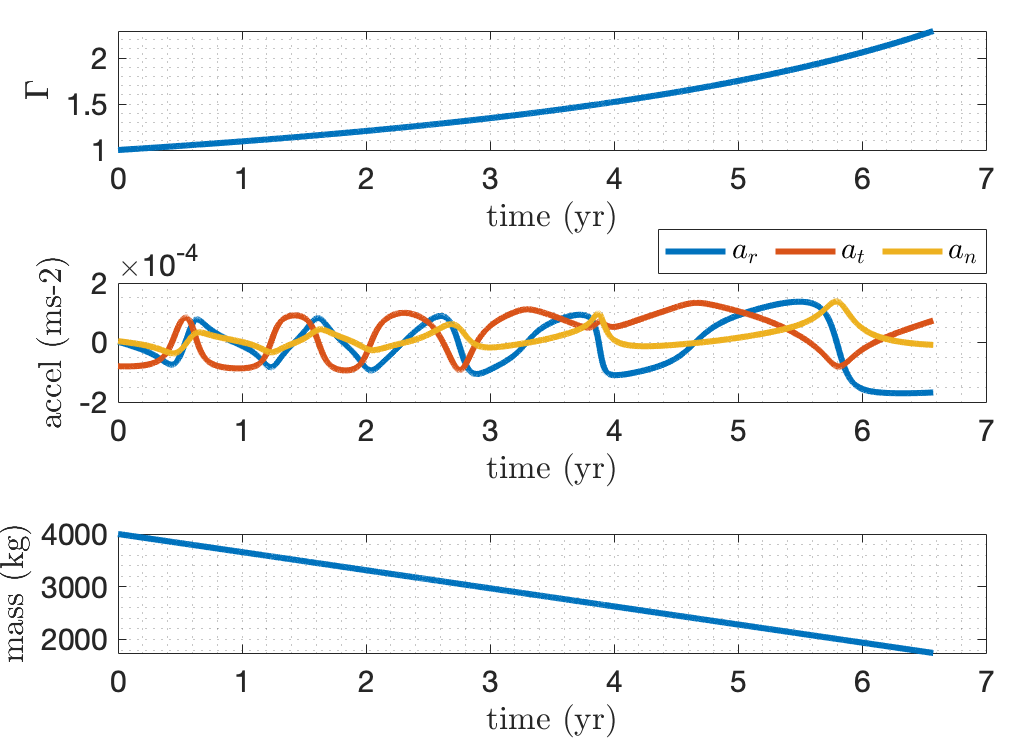} }}%
    \caption{Time-optimal Earth to Dionysus transfer}%
  
\end{figure}

 {Tables \ref{topt2} and \ref{tsol} show the TO results obtained for the Earth to Dionysus test case.} Figure \ref{dta} shows the time-optimal trajectory while Fig. \ref{dtb} shows the profiles of thrust acceleration and spacecraft mass along the trajectory. 

\begin{table}[!hbt]
\centering 
\caption{Initial costates obtained at each step of the time-optimal solution for the Earth to Dionysus trajectory}
\begin{tabular}{llccccccc}
\toprule
Method                                                           & Dynamic& TOF                      & $\lambda_p$                       & $\lambda_f$                                                             & $\lambda_g$                         & $\lambda_h$                         & $\lambda_k$                         & $\lambda_L$                         \\
~  & Model              &days & ${L_u/T_u^3}$ & ${L_u^2/T_u^3}$& ${L_u^2/T_u^3}$ & ${L_u^2/T_u^3}$  & ${L_u^2/T_u^3}$& ${L_u^2/T_u^3}$ \\ \hline
EO & Kep. &   3534.0000    &-1.7649  &  -0.2215   &  1.0965 &   -1.0684 &   -2.3545 &   -0.0096 \\
EO (for TOF guess) & Kep. & 2098.0432    &  4.7057  & -0.2692 &   3.7786  & -3.7088  & -4.1064 &  -0.1425\\
TO  & Kep. &   2401.4301  & 2.1022  & -0.3377  &  1.8991  & -1.4481  & -0.9391&   -0.1133 \\ \bottomrule
\end{tabular}
\label{topt2}
\end{table}


\subsection{Transfer between Space Debris}
The framework for this test case comes from the 9\textsuperscript{th} edition of the Global Trajectory Optimisation Competition (GTOC), centered around space debris removal. In \cite{Armellin2018c} test cases from GTOC 9 are recreated assuming  {low-thrust} propulsion instead of impulsive thrust. Here, case F given in \cite{Armellin2018c} is solved where a transfer between two synthetic space debris occurs.  {The boundary states for this transfer are $\boldsymbol{x}_0 = [1.117658, -0.000418, 0.000555, -1.040879, -0.511994, 1.706348]^T$ and $\boldsymbol{x}_1 = [1.123581, 0.001719, 0.007079, -1.025061, -0.525796, 87.229928]^T$ in MEE coordinates.} 
 As this is a geocentric case, the effect of $J_2$ and eclipses are also considered. Similarly to \cite{singh2021eclipse}, the effect of eclipses is accounted for by switching the propulsion system off without deriving the midpoint boundary conditions.   
\subsubsection{Fuel-Optimal Solution}
The effect of $J_2$ is introduced at the start of the problem when solving the FO problem. Then, a secondary continuation method is used to introduce eclipses. This method is given in Algorithm \ref{eclipses}. 
\begin{algorithm}[hbt!]
\caption{Introducing Eclipses and $J_2$}\label{eclipses}
\begin{algorithmic}
\State Solve the FO problem with discontinuous thrust and the influence of $J_2$ using Algorithm \ref{method} to obtain $\boldsymbol{\lambda}_{f_{J_2}}$
\State Let $\epsilon = 0$ and $\boldsymbol{\lambda}_{g} = \boldsymbol{\lambda}_{f_{J_2}}$
\While{$\epsilon \leq 1$}
    \State In FO dynamics, let $\Gamma_f = \Gamma_f (1 - \epsilon \nu)$ \Comment{$\nu$ is determined as described in \cite{Jonathan2019},  with $c_t = 100, c_s = 0.9$.}
    \State Solve the FO problem using $\boldsymbol{\lambda}_{g}$ as an initial guess to get solution $\boldsymbol{\lambda}_{f_{J_2+ecl}}$.
    \State  $\boldsymbol{\lambda}_{g} = \boldsymbol{\lambda}_{f_{J_2+ecl}}$ and  $\epsilon = \epsilon + \Delta \epsilon$
    \EndWhile

\end{algorithmic}
\end{algorithm}
Note that the largest value of $ \Delta \epsilon$ that allowed convergence at all continuation steps was 0.1.\par
\begin{figure}[!hbt]
    \centering
    \subfloat[\centering Optimal Trajectory    \label{geoplots1a}]{{\includegraphics[width=8.0cm]{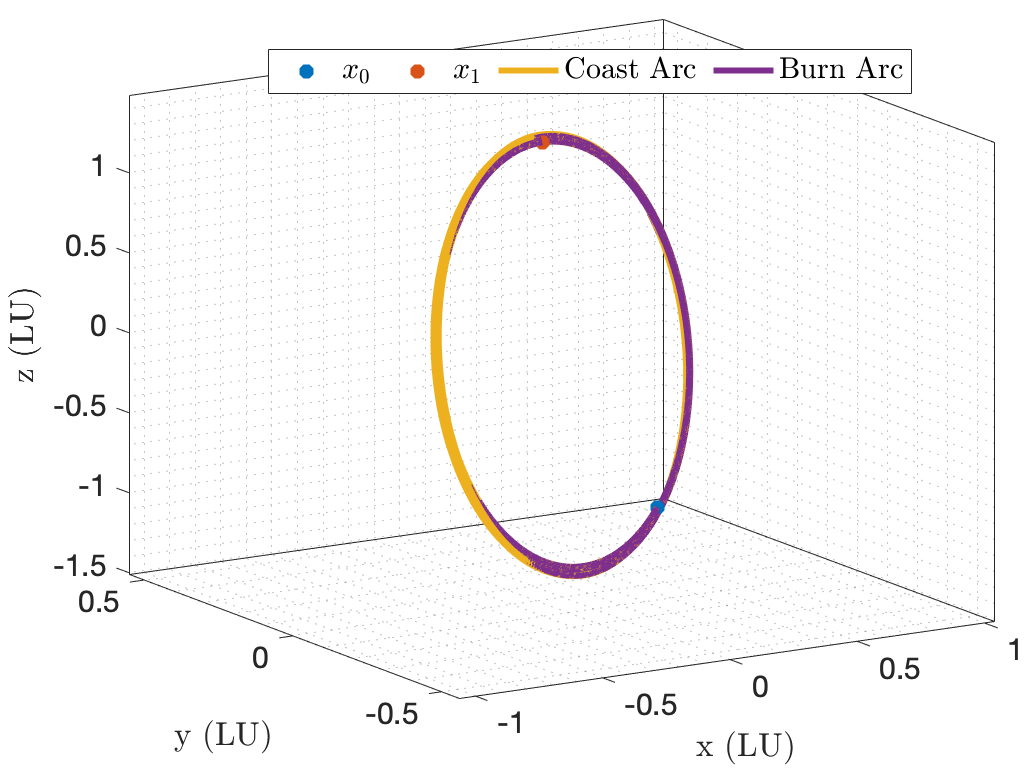} }}%
    \qquad
    \subfloat[\centering Control and spacecraft mass profiles \label{geoplots1b}]{{\includegraphics[width=7.5cm]{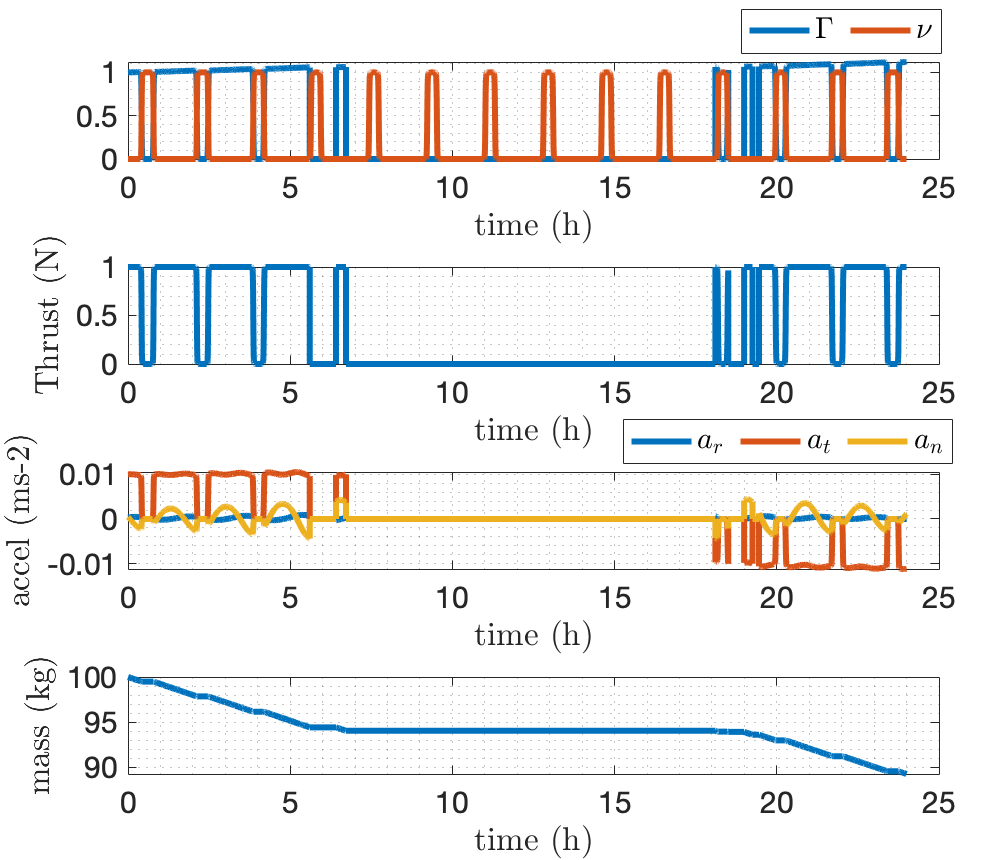} }}%
    \caption{Fuel-optimal debris to debris transfer}%
 
\end{figure}

 
 {Tables \ref{geot1} and \ref{tsol} show the FO results obtained, for a 1 day time of flight, starting on 31 Dec. 2023, 00:00:00.0 UTC}. Five continuation steps are required to generate the FO solution under the effect of $J_2$ and Keplerian dynamics. Then ten continuation steps are taken to introduce the full effect of eclipses. Figure \ref{geoplots1a} shows the FO trajectory, which has multiple revolutions and many coast and burn arcs. It can be seen that the coast arcs are mainly accumulated towards the left side of Fig. \ref{geoplots1a} due to eclipses. Figure \ref{geoplots1b} illustrates the profiles of $\Gamma$ and the spacecraft's mass. When the Keplerian dynamics and $J_2$ are considered, the optimal $\Delta v$ obtained is \SI{317.58}{\meter\per\second}, which agrees with the value reported in \cite{Armellin2018c}.

\begin{table}[!hbt]
\centering 
\caption{Initial costates obtained at each step of the fuel optimal solution for the GTOC 9 transfer}
\begin{tabular}{llccccccc}
\toprule
Method & Dynamic & Fuel mass & $\lambda_p$ & $\lambda_f$& $\lambda_g$ & $\lambda_h$  & $\lambda_k$& $\lambda_L$\\
~  & Model              &kg & ${L_u/T_u^3}$ & ${L_u^2/T_u^3}$& ${L_u^2/T_u^3}$ & ${L_u^2/T_u^3}$  & ${L_u^2/T_u^3}$& ${L_u^2/T_u^3}$ \\
\hline
EO & Kep. + $J_2$ & 12.5444  & -41.3247&-2.6164 &-1.3320 &-10.2564 & -1.4592 &0.6874 \\
FO & Kep. + $J_2$ & 10.2328  & -33.7350& -2.0584 & -1.2709&-15.2377 &-2.2941 & 0.5531         \\
FO & Kep. + $J_2 $ $+$ eclipses  & 11.0908    & -40.7995 &-2.9911& -2.4330&-16.6356 & -5.5040 &0.6692  \\ \bottomrule
\end{tabular}
\label{geot1}
\end{table}

\subsubsection{Time-Optimal Solution} 
\begin{figure}[!hbt]
    \centering
    \subfloat[\centering Optimal Trajectory  \label{geotimea}]{{\includegraphics[width=7.5cm]{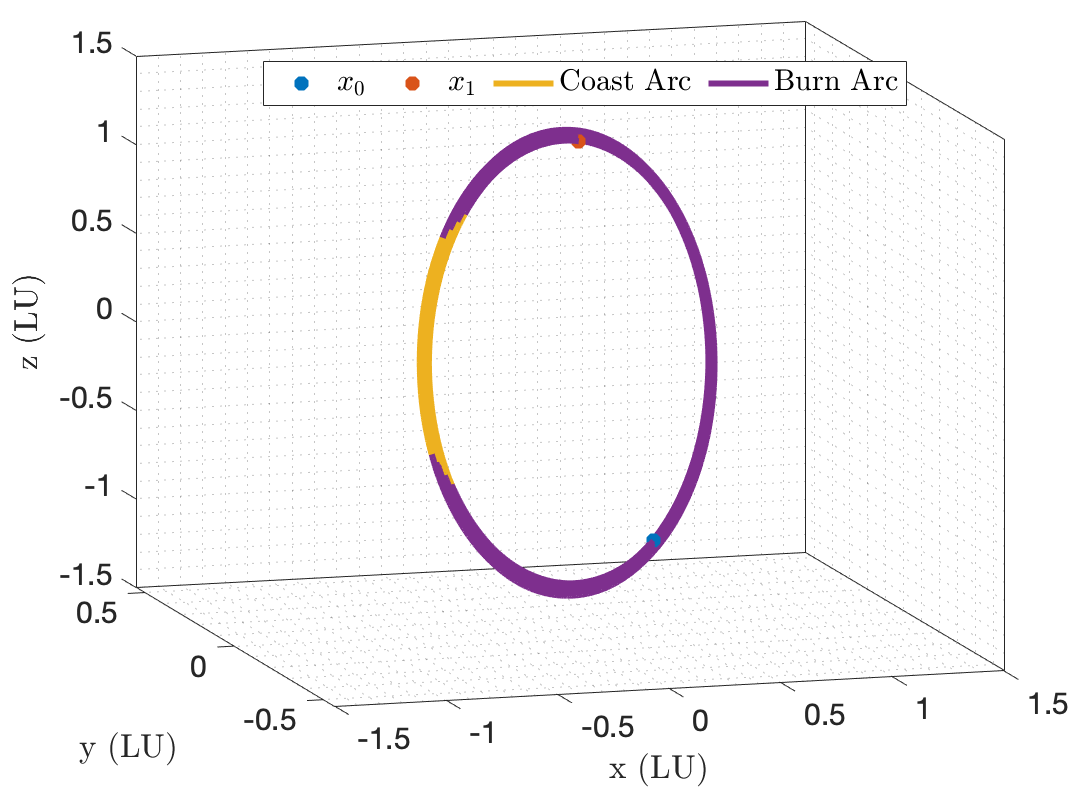} }}%
    \qquad
    \subfloat[\centering Control and spacecraft mass profiles  \label{geotimeb}]{{\includegraphics[width=7.5cm]{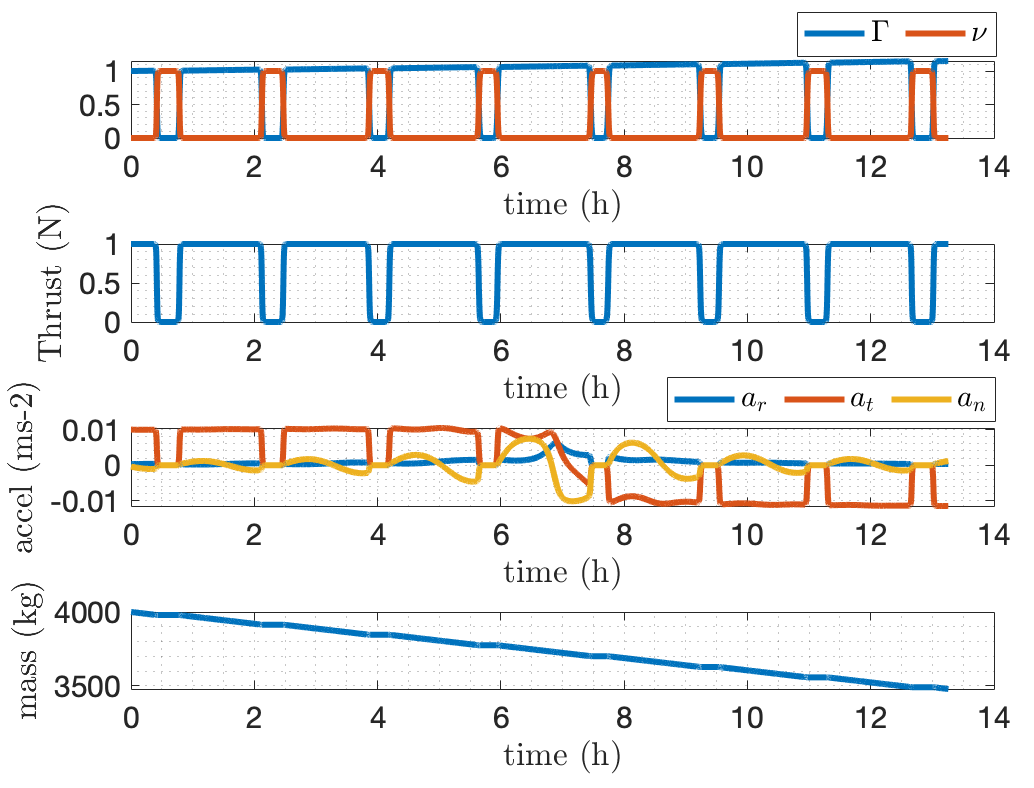} }}%
    \caption{Time-optimal debris to debris transfer}%
\end{figure}

When solving the TO problem, it was noted that introducing $J_2$ as the final step resulted in a higher convergence rate. Hence, the eclipses were introduced first, followed by $J_2$.  {Tables \ref{toptgeo1} and \ref{tsol} show the TO results obtained}. $\beta_t$ was recalculated at each iteration of the eclipse introduction to improve the convergence speed further. Figure \ref{geotimea} shows the time-optimal trajectory while
Fig. \ref{geotimeb} illustrates the profiles of thrust acceleration and spacecraft mass. 

\begin{table}[!hbt]
\centering 
\caption{Initial costates obtained at each step of the time-optimal solution for the GTOC9 trajectory}
\begin{tabular}{llccccccc}
\toprule
Method                                                           & Dynamic& TOF                      & $\lambda_p$                       & $\lambda_f$                                                             & $\lambda_g$                         & $\lambda_h$                         & $\lambda_k$                         & $\lambda_L$                         \\
~  & Model              &days & ${L_u/T_u^3}$ & ${L_u^2/T_u^3}$& ${L_u^2/T_u^3}$ & ${L_u^2/T_u^3}$  & ${L_u^2/T_u^3}$& ${L_u^2/T_u^3}$ \\
\hline
EO & Kep. &     1.0000   &   -41.3247&-2.6164 &-1.3320 &-10.2564 & -1.4592 &0.6874        \\
EO & Kep. (for TOF guess)  &   0.5625 &   -85.5442  &  -6.9855 &   -2.7904 &  -35.2564 &   32.1044 &    2.5484\\
TO & Kep. &  0.5572 & -77.9743 &   -5.2764   & -2.8717 &   -21.0868 &    19.7377   &  2.2013    \\
TO & Kep.$+$ eclipses  &0.5499& -70.0088   & -5.0282 &   -4.2352 &  -22.3886 &   20.4769 &    1.9959\\
TO & Kep. + $J_2$ + eclipses   & 0.5525 &  -141.9498  &  -9.7304 &   -3.7754 &  -28.0852 &   10.6069  &   4.0443 \\ \bottomrule
\end{tabular}
\label{toptgeo1}
\end{table}



\begin{table}[hbt!]

\caption{ {Time- and fuel-optimal solutions}}

\begin{tabular}{l|cccc|cccc}
\toprule
~ & \multicolumn{4}{c|}{Fuel Optimal}                                                                                                                                                                 & \multicolumn{4}{c}{Time Optimal}                                                                                                                                                                  \\\hline 
 Mission    & \begin{tabular}[c]{@{}c@{}}Flight\\ Time (d)\end{tabular} & \begin{tabular}[c]{@{}c@{}}Fuel Cons. \\ (kg)\end{tabular} & \begin{tabular}[c]{@{}c@{}}Comp. \\  Time (s)\end{tabular} & $\Gamma_{TR}$ & \begin{tabular}[c]{@{}c@{}}Flight \\ Time (d)\end{tabular} & \begin{tabular}[c]{@{}c@{}}Fuel Cons. \\ (kg)\end{tabular} & \begin{tabular}[c]{@{}c@{}}Comp. \\  Time (s)\end{tabular} &$\beta_t$ \\\hline 
Earth to Tempel 1     & 420.00                                                       & 348.26                                                     & 0.114                                                     & 0.4781  & 327.15                                                        & 576.50                                                     & 0.208                                                     & 19.9859 \\
Earth to Dionysus    & 3534.00                                                      & 1280.70                                                    & 0.425                                                     & 0.5389  & 2401.43                                                       & 2256.80                                                    & 0.733                                                      & 1.5928  \\
Space debris    & 1.00                                                         & 11.09                                                      & 75.600  & 0.4886  & 0.553                                                         & 13.03                                                      & 55.320 & 20.9637 \\
\bottomrule
\end{tabular}
\label{tsol}
\end{table}

Overall, it can be seen that the discussed methodology can solve heliocentric transfers and geocentric transfers with eclipses and perturbations to a high degree of accuracy within a short amount of time. Implementing the $\Gamma_{TR}$ constant in the fuel optimization and $\beta_t$ in the time optimization ensures convergence and reduces the number of iterations required compared to the classical case where both these parameters are set to unity. By integrating mass separately in both fuel and time-optimal cases, the computational speed is shown to improve. 

\section{Conclusion}
This note has proposed techniques to improve the convergence of indirect trajectory optimization methods when solving fuel- and time-optimal problems. By introducing a constant ($\Gamma_{TR}$) in the fuel optimal objective function, a good transition from energy-optimal (EO) to fuel-optimal (FO) problem is achieved. This approach avoided the need for random initial guesses when using control smoothing and allowed a low number of iterations to reach convergence. Similarly, it was shown that EO could provide a good guess for TO problems and that introducing the constant $\beta_t$ to the objective function can reduce the initial residuals and improve convergence. Compared with the traditional continuation methods where $\Gamma_{TR}$ and $\beta_t$ are taken to be unity, our approach brings advantages in both the computational efficiency and convergence rate.\par 
Furthermore, it was shown that the equation for the Lagrangian multiplier of the mass and the associated boundary condition could be ignored for both the FO and TO cases without affecting optimality. This simplification reduces the problem dimension and thus improves efficiency.
Numerical examples are given to validate the advantages and the robustness of the discussed techniques for heliocentric and geocentric cases. For geocentric cases, the effect of eclipses and $J_2$ perturbations were also considered.

\bibliographystyle{unsrt}

\begin{thebibliography}{99}

\bibitem{1}
Xun Pan and Binfeng Pan.
\newblock Practical homotopy methods for finding the best minimum-fuel transfer
  in the circular restricted three-body problem.
\newblock {\em IEEE Access journal}, PP:1--1, 03 2020.

\bibitem{2}
Anastassios Petropoulos and Jon Simst.
\newblock A review of some exact solutions to the planar equations of motion of
  a thrusting spacecraft.
\newblock In {\em Proceedings of the 2nd International Symposium Low Thrust
  Trajectories, Toulouse, France}, pages 10--30, 2002.

\bibitem{3}
David Morante, Manuel Sanjurjo~Rivo, and Manuel Soler.
\newblock A survey on low-thrust trajectory optimization approaches.
\newblock {\em MDPI journal of Aerospace}, 8(3), 2021.

\bibitem{Geller2017}
David~K Geller, Nicholas Ortolano, and Aaron Avery.
\newblock Autonomous optimal trajectory planning for orbital rendezvous,
  satellite inspection, and final approach based on convex optimization.
\newblock {\em The Journal of the Astronautical Sciences}, 68:444--479, 2021.

\bibitem{Taheri2021}
Ehsan Taheri.
\newblock Optimization of many-revolution minimum-time low-thrust trajectories
  using sundman transformation.
\newblock {\em AIAA Scitech 2021 Forum}, pages 1--17, 2021.

\bibitem{Wu20212}
Di~Wu, Wei Wang, Fanghua Jiang, and Junfeng Li.
\newblock Minimum-time low-thrust many-revolution geocentric trajectories with
  analytical costates initialization.
\newblock {\em Aerospace Science and Technology journal}, 119:107--146, 12
  2021.

\bibitem{Altitude-dependent2010c}
Di~Wu, Fanghua Jiang, and Junfeng Li.
\newblock Warm start for low-thrust trajectory optimization via switched
  system.
\newblock {\em Journal of Guidance, Control, and Dynamics}, 47:8--11, 2021.

\bibitem{Wu2021}
Di~Wu, Lin Cheng, Fanghua Jiang, Junfeng Li, and L~Cheng.
\newblock Rapid generation of low-thrust many revolution earth-center
  trajectories based on analytical state-based control.
\newblock {\em Acta Astronautica}, 2021.

\bibitem{Rasotto2015}
M~Rasotto, R~Armellin, and P~Di Lizia.
\newblock Multi-step optimization strategy for fuel-optimal orbital transfer of
  low-thrust spacecraft.
\newblock {\em Acta Astronautica}, 13:1--27, 2015.

\bibitem{Kechichian1997}
Jean~Albert Kechichian.
\newblock Reformulation of edelbaum's low-thrust transfer problem using optimal
  control theory.
\newblock {\em Journal of Guidance, Control, and Dynamics}, 20:988--994, 1997.

\bibitem{Woollands2020}
Robyn Woollands, Ehsan Taheri, and John~L. Junkins.
\newblock Efficient computation of optimal low thrust gravity perturbed orbit
  transfers.
\newblock {\em Journal of the Astronautical Sciences}, 67:458--484, 6 2020.

\bibitem{Jin2020}
Kai Jin, David~K. Geller, and Jianjun Luo.
\newblock Robust trajectory design for rendezvous and proximity operations with
  uncertainties.
\newblock {\em Journal of Guidance, Control, and Dynamics}, 43:741--753, 2020.

\bibitem{Ravikumar2020}
L.~Ravikumar, Radhakanth Padhi, and N.K. Philip.
\newblock Trajectory optimization for rendezvous and docking using nonlinear
  model predictive control.
\newblock {\em IFAC-PapersOnLine}, 53(1):518--523, 2020.

\bibitem{Armellin2018c}
Roberto Armellin, David Gondelach, Juan~Félix San-Juan, and Juan Felix~San
  Juan.
\newblock Multiple revolution perturbed lambert problem solvers.
\newblock {\em Journal of Guidance, Control, and Dynamics}, 41:2019--2032, 6
  2018.

\bibitem{Russell2007a}
Ryan~P. Russell.
\newblock Primer vector theory applied to global low-thrust trade studies.
\newblock {\em Journal of Guidance, Control, and Dynamics}, 30:460--472, 2007.

\bibitem{Wall2009}
Bradley~J. Wall and Bruce~A. Conway.
\newblock Shape-based approach to low-thrust rendezvous trajectory design.
\newblock {\em Journal of Guidance, Control, and Dynamics}, 32:95--102, 2009.

\bibitem{prussing_2010}
John~E. Prussing.
\newblock {\em Primer Vector Theory and Applications}.
\newblock Cambridge University Press, 2010.

\bibitem{Kim2008}
Seung-Jean Kim, K.~Koh, M.~Lustig, Stephen Boyd, and Dimitry Gorinevsky.
\newblock An interior-point method for large-scale l1-regularized least
  squares.
\newblock {\em Selected Topics in Signal Processing, IEEE Journal of}, 1:606 --
  617, 01 2008.

\bibitem{Taheri2018}
Ehsan Taheri and John~L. Junkins.
\newblock Generic smoothing for optimal bang-off-bang spacecraft maneuvers.
\newblock {\em Journal of Guidance, Control, and Dynamics}, 41(11):2470--2475,
  2018.

\bibitem{Jiang2012}
Fanghua Jiang.
\newblock Practical techniques for low-thrust trajectory optimization with
  homotopic approach.
\newblock {\em Journal of Guidance Control and Dynamics}, 1:245--258, 01 2012.

\bibitem{Bertrand}
R.~Bertrand and R.~Epenoy.
\newblock New smoothing techniques for solving bang–bang optimal control
  problems—numerical results and statistical interpretation.
\newblock {\em Optimal Control Applications and Methods}, 23(4):171--197, 2002.

\bibitem{Haberkorn}
T.~Haberkorn, P.~Martinon, and J.~Gergaud.
\newblock Low thrust minimum-fuel orbital transfer: A homotopic approach.
\newblock {\em Journal of Guidance, Control, and Dynamics}, 27(6):1046--1060,
  2004.

\bibitem{mansell2018adaptive}
Justin~R Mansell and Michael~J Grant.
\newblock Adaptive continuation strategy for indirect hypersonic trajectory
  optimization.
\newblock {\em Journal of Spacecraft and Rockets}, 55(4):818--828, 2018.

\bibitem{Martinon2007}
Pierre Martinon and Joseph Gergaud.
\newblock Using switching detection and variational equations for the shooting
  method.
\newblock {\em Optimal Control Applications and Methods}, 28, 03 2007.

\bibitem{tang2018fuel}
Gao Tang, Fanghua Jiang, and Junfeng Li.
\newblock Fuel-optimal low-thrust trajectory optimization using indirect method
  and successive convex programming.
\newblock {\em IEEE Transactions on Aerospace and Electronic Systems},
  54(4):2053--2066, 2018.

\bibitem{ZHU201798}
Zhengfan Zhu, Qingbo Gan, Xin Yang, and Yang Gao.
\newblock Solving fuel-optimal low-thrust orbital transfers with bang-bang
  control using a novel continuation technique.
\newblock {\em Acta Astronautica}, 137:98--113, 2017.

\bibitem{Jiang2017}
Fanghua Jiang, Gao Tang, and Junfeng Li.
\newblock Improving low-thrust trajectory optimization by adjoint estimation
  with shape-based path.
\newblock {\em Journal of Guidance, Control, and Dynamics}, 40:3280--3287,
  2017.

\bibitem{Ayyanathan2022}
Praveen~J. Ayyanathan and Ehsan Taheri.
\newblock Mapped adjoint control transformation method for low-thrust
  trajectory design.
\newblock {\em Acta Astronautica}, 4 2022.

\bibitem{palermoCAM}
Maria~Francesca Palermo, Andrea De~Vittori, Pierluigi Di~Lizia, and Roberto
  Armellin.
\newblock Low-thrust collision avoidance maneuver optimization.
\newblock {\em Journal of Guidance, Control, and Dynamics}, 2022.
\newblock (accepted for publication).

\bibitem{Betts2000VeryLT}
John~T. Betts.
\newblock Very low-thrust trajectory optimization using a direct sqp method.
\newblock {\em Journal of Computational and Applied Mathematics}, 120:27--40,
  2000.

\bibitem{don}
Eelco Doornbos.
\newblock {\em Thermospheric Density and Wind Determination from Satellite
  Dynamics}.
\newblock Springer, 01 2012.

\bibitem{Vallado2013}
David~A Vallado and Wayne~D Mcclain.
\newblock {\em Fundamentals of Astrodynamics and Applications Fourth Edition}.
\newblock Space Technology Library, 2013.

\bibitem{Jonathan2019}
Jonathan Aziz, Daniel Scheers, Jeffrey Parker, and Jacob Englander.
\newblock A smoothed eclipse model for solar electric propulsion trajectory
  optimization.
\newblock {\em Transactions of the Japan Society for Aeronautical and Space
  Sciences, Aersoapce Technology Japan}, 17(2):181--188, 2019.

\bibitem{Lizia2014}
Pierluigi Di~Lizia, Roberto Armellin, Alessandro Morselli, and Franco
  Bernelli-Zazzera.
\newblock High order optimal feedback control of space trajectories with
  bounded control.
\newblock {\em Acta Astronautica}, 94, 02 2014.

\bibitem{galassi2002gnu}
Brian Gough.
\newblock {\em GNU Scientific Library Reference Manual}.
\newblock Network Theory Ltd., 3rd edition, 2009.

\bibitem{nasa_2021}
Gridded ion thrusters (next-c) - glenn research center, Dec 2021.

\bibitem{satsearch}
B20 thruster - green chemical propulsion, 2022.

\bibitem{singh2021eclipse}
Sandeep Singh, John Junkins, Brian Anderson, and Ehsan Taheri.
\newblock Eclipse-conscious transfer to lunar gateway using ephemeris-driven
  terminal coast arcs.
\newblock {\em Journal of Guidance, Control, and Dynamics}, 44(11):1972--1988,
  2021.

\end{thebibliography}

\end{document}